\documentclass{gen-p-l}
\usepackage{amsfonts}
\usepackage{amssymb}
\usepackage{amsthm}
\usepackage{upref}
\makeatletter
\def\LaTeX{\leavevmode L\raise.42ex
    \hbox{\kern-.3em\size{\sf@size}{0pt}\selectfont A}\kern-.15em\TeX}
 \makeatother
 
\numberwithin{equation}{section}

\newtheorem{lemma}{Lemma}[section]
\newtheorem{theorem}[lemma]{Theorem} 
\newtheorem{corollary}[lemma]{Corollary}
\newtheorem{proposition}[lemma]{Proposition}

\theoremstyle{definition}

\newtheorem{remark}[lemma]{Remark}

\renewcommand{\det}{\operatorname{Det}}
\newcommand{\tr}{\operatorname{Tr\,}}

  \newcommand{\slim}{\operatorname{s-lim}}

\newcommand{\HH}{\mathsf{H}}
 
  \newcommand{\e}{\eqref}
\newcommand{\ri}{\rightarrow}
\newcommand{\q}{\quad}

   \newcommand{\sgn}{\operatorname{sgn}}
  \newcommand{\ran}{\operatorname{Ran}}

       \newcommand{\var}{\operatorname{var}}

\renewcommand\Im{\operatorname{Im}}
\renewcommand\Re{\operatorname{Re}}

\newenvironment{pf}{\begin{proof}}{\end{proof}}

\def\qqq{\mathrel{\subset\mkern-15mu\lower.38ex\hbox{${\scriptscriptstyle\rightarrow}$}}}

\let\cal\mathcal

\let\Bbb\mathbb

\begin{document}
\title[fourth order differential operators]
{Spectral and scattering theory \\ of fourth order differential operators}
\author{ D. R. Yafaev}
\address{ IRMAR, Universit\'{e} de Rennes I\\ Campus de
  Beaulieu, 35042 Rennes Cedex, FRANCE}
\email{yafaev@univ-rennes1.fr}

\dedicatory{To Mikhail Shl\"emovich Birman on his 80-th birthday}

\begin{abstract}
  An ordinary differential operator
of the fourth order   with coefficients converging at infinity sufficiently rapidly  to constant limits  is considered. Scattering theory for this operator is developed in terms of special solutions of the corresponding differential  equation. In contrast to equations of second order ``scattering" solutions contain exponentially decaying terms. A relation between the scattering matrix and a matrix of coefficients at exponentially decaying modes is found.
In the second part of the paper   the operator $D^4$ on the half-axis with different  boundary conditions at the point zero is studied. Explicit formulas for basic objects of the scattering theory are found. In particular, a classification of different types of zero-energy resonances is given.
\end{abstract}

\maketitle

\thispagestyle{empty}

\section{Introduction}  
{\bf 1.1. }
General scattering theory for differential operators does not depend on the order of  
  operators (see \cite{B} for the trace class approach and \cite{Ku} for the  smooth approach). Suppose that coefficients of a differential operator $H$ converge sufficiently rapidly to constant values at infinity, and let $H_{0}$ be the operator with these constant coefficients. Then the wave operators $W_{\pm}(H,H_{0})$ for the pair $H_{0}$, $H$ exist  and are complete, and  the corresponding scattering matrix is a unitary operator.
 The operator $H$ does not have the singular continuous spectrum, and its point spectrum might accumulate only at thresholds (critical values of the symbol of the operator $H_{0}$). Moreover, an expansion theorem in eigenfunctions of the operator $H$ is true.
 
  However the behavior of     eigenfunctions 
of  the  continuous spectrum at infinity is essentially simpler
  for differential operators of the second order   than for higher order differential operators. This is intimately related to the fact that higher order differential operators might have eigenvalues embedded in the continuous spectrum.
  
\medskip

{\bf 1.2. } 
 In this paper we consider one-dimensional differential  operators. For definiteness, we choose operators of the fourth order.
  Set $H_{0}=D^4$ and
\begin{equation}
H=D^4+ D v_{1} (x) D +v_{0}(x),\q D= -i d /dx,
 \label{eq:2.1a}\end{equation}
 where  the  functions $v_{0}$ and $v_{1}$ are real and satisfy a short-range  
 assumption  
 \begin{equation}
 v_{j}\in L_1({\Bbb R}),\q j=0,1.
\label{eq:2.1}\end{equation}
A smoothness of the function $v_{1}(x)$ is not required.
 We are interested in construction of the wave operators and
of the scattering matrix in terms of solutions of the corresponding differential equation
\begin{equation}
( u '''(x)-  v_{1}(x)u '(x))' + v_{0}(x)u  (x) = \lambda u(x), \q \lambda>0. 
\label{eq:dex}\end{equation}
 This construction turns out to be  more complicated than for equations of the second  order. Indeed,  the ``free" equation $u^{(4)}(x) =\lambda  u (x)$    has solutions $e^{ik x}$, $e^{-ik x}$, $e^{-k x}$ and $e^{k x}$ where $\lambda=k^4$, $k>0$. If, for example, the coefficients $v_{0}(x)$ and $v_1(x)$ are compactly supported, then every solution of equation \e{eq:dex} is a linear combination of these exponentials for large positive and large negative $x$.   Eigenfunctions $\psi_{j} (x, \lambda)$, $j=1,2 $, 
of  the  operator $H$ (its continuous spectrum has multiplicity two) are special solutions of equation \e{eq:dex}.
It is natural to expect that   they 
do  not contain exponentially
increasing terms  and that the scattering matrix is determined only
by coefficients at oscillating terms $e^{ik x}$ and
$e^{-ik x}$. We justify this conjecture   under general short-range assumption \e{eq:2.1}.

Furthermore, if the coefficients $v_{1} (x)$ and $v_{0}(x)$ decay  super-exponentially at infinity, then it is possible to distinguish exponentially
decreasing terms  in the asymptotics of the functions $\psi_{j} (x, \lambda)$, $j=1,2 $,  as $x\to \pm\infty$. Our main observation in the first part of the paper is that the coefficients at these terms determine the scattering matrix. This might be eventually of interest for a study of the (inverse) problem of a reconstruction of the coefficients $v_{0} (x)$ and $v_{1}(x)$ from scattering data. We refer to \cite{BDT} for a comprehensive study of the inverse problem. 

Recall that for differential operators of the second order
the eigenfunctions of the continuous spectrum can be constructed (see \cite{F})   with a help of Volterra integral equations. This procedure seems not to work  for 
operators of higher order. Therefore we use a general scheme
of scattering theory (see, e.g.,  \cite{Kuroda, Ya}) relying on the Lippmann-Schwinger equation.

In Section~2 we present a stationary approach to scattering theory for the operator $H$. This approach is quite general and, up to some technical details, works for     multi-dimensional differential operators of an arbitrary order (see \cite{Ku}). The absence of the singular continuous spectrum is verified in Section~3. Here  instead of Agmon's bootstrap arguments \cite{Ag}, a method specific for ordinary differential equations is used.  
 After this prerequisite, we study in Section~4 asymptotic behavior of eigenfunctions of the continuous spectrum. The relation between the scattering matrix and the coefficients at   exponentially decreasing modes is formulated in Theorem~\ref{4di3}.

\medskip

{\bf 1.3. } 
The second part of the paper (Sections~5, 6  and 7) is devoted to a study of the operator $H=D^4 $  in the space $L_{2}({\Bbb R}_{+})$
with some self-adjoint boundary conditions at the point $x=0$.  This operator  can be compared with a well known Hamiltonian $ D^2$ with a boundary condition 
\begin{equation}
u'(0) =\alpha u(0), \q \alpha=\bar{\alpha}. 
\label{eq:BCT}\end{equation}
 This boundary condition is interpreted as a point interaction at the point $x=0$. The point interaction is a good approximation to a pertutbation by  a potential for low (but not for high) energies.

  We shall  write down explicit (although not very simple) formulas for basic objects of  
scattering theory for the operator $H$, such as the resolvent kernel, eigenfunctions, the scattering matrix, the perturbation determinant, the spectral shift function and so forth. This model seems to be of interest because it allows us to analyse a behavior of different objects at zero energy as well as at a positive eigenvalue of
the operator $H$.  In particular, we discuss different types of zero-energy resonances.

 There exist several definitions of zero-energy resonances which are essentially (but not completely) equivalent.   According to a general variational  definition
  (see \cite{Bia}) an operator $H$ has a zero-energy resonance if for a small negative perturbation an additional negative eigenvalue appears. Of course this definition depends not only on $H$ but also on a class of perturbations.
   Other definitions (see \cite{Y15} or \cite{LNM}) are adapted to differential operators. Thus, 
    the Schr\"odinger operator $H= -\Delta + v(x)$ has a zero-energy resonance if the kernel of its resolvent $R(z)= (H-z)^{-1}$ has a singularity  $\varphi (x) \overline{\varphi(x')}(-z)^{-1/2}$ as $|z| \to 0$.  Here $\varphi (x)$ is a solution of the equation $H \varphi=0$ which is bounded at infinity (in the one-dimensional case). The existence of such solutions gives still another criterium for  the appearance of zero-energy resonances. Zero-energy resonances can be considered as a weakened version of bound states (for zero energy) and are often called half-bound states. This point of view is confirmed by the behavior of the spectral shift function  $\xi(\lambda)$  
        at the point $\lambda =0$. Consider, for example, the operators   $H_{0} $ and $ H(\alpha)$ corresponding to the differential expression $D^2$ in the space $L_{2}({\Bbb R}_{+})$ with boundary conditions $u(0)=0$ and \e{eq:BCT}, respectively. The operator $H(\alpha)$ has a zero-energy resonance if and only if $\alpha=0$. The spectral shift function  $\xi(\lambda)$ for the pair $H_{0} $, $ H(\alpha)$ is continuous at $\lambda=0$ if $\alpha\neq 0$  and it has the  jump   $-1/2$  if $\alpha=0$. This should be compared with the fact that  the jump of  $\xi(\lambda)$ equals $-1$
        at an   eigenvalue of   $H$.

In Section~7  we       analyse in some details     zero-energy resonances for the operator $H = D^4 $ with different boundary conditions at the point $x=0$.
  It turns out that, compared to second order, for fourth order differential operators their notion acquires some new featutes although the variational definition remains of course  valid. A more detailed classification is given in terms of singularities of the resolvent at the point $z=0$, of a behavior of solutions (which are of course all polynomials of degree $3$) of the equation $u^{(4)}(x)=0$ satisfying the boundary condition and of a jump of the corresponding spectral shift function at the point $\lambda =0$. 
   This analysis shows that it is natural to introduce now $1/4$- and $3/4$-bound states. Here we mention only that the operator $H$ has a $1/4$-bound state (a $3/4$-bound state) if a non-trivial linear function (a constant) satisfies   the boundary conditions at $x=0$.
   
    \medskip

Thanks are due to V. Suhanov for a useful discussion.


\section{General scattering theory }

{\bf 2.1. }
Under the assumption 
\[
\sup_{x\in {\Bbb R}}\int_x^{x+1} (|v_{1} (y)|  + |v_{0} (y)|  )dy <\infty
\]
 we have that, for all $\varepsilon>0$ and all functions $u$
from the Sobolev class ${\HH}^2 ({\Bbb R})$,
\begin{eqnarray}
\int_{-\infty}^\infty (|v_{1} (x)| |u'(x)|^2+ |v_{0} (x)| |u(x)|^2)dx
\nonumber\\
\leq
\varepsilon \int_{-\infty}^\infty |u''(x)|^2 dx + C(\varepsilon)  \int_{-\infty}^\infty |u(x)|^2 dx.
\label{eq:QF}\end{eqnarray} 
Here and below $C$ denotes different positive constants whose precise values are of no importance. Therefore the   quadratic form
\[
h[u,u] = \int_{-\infty}^\infty (|u''(x)|^2+ v_{1} (x)|u'(x)|^2+v_{0} (x)|u(x)|^2)dx
\]
is semibounded from below and is closed on   ${\HH}^2 ({\Bbb R})$.
Thus, it defines a self-adjoint    operator $H$    in the Hilbert space ${\cal H}=L_2({\Bbb R})$ with domain ${\cal D}(H)\subset {\HH}^2 ({\Bbb R})$. This operator   corresponds  to  formal differential expression \e{eq:2.1a}.

Note that the operators  $H_0=D^4 : {\HH}^2 ({\Bbb R})\to {\HH}^{-2} ({\Bbb R})$ and, by virtue of estimate \e{eq:QF},
\begin{equation}
V= D v_{1}(x) D +v_{0}(x) : {\HH}^2 ({\Bbb R})\to {\HH}^{-2} ({\Bbb R})
\label{eq:V}\end{equation}
are bounded operators. It can easily be shown that a function $u\in {\HH}^2 ({\Bbb R})$ belongs to ${\cal D}(H)$ if and only if $H_{0}u+ Vu \in L_{2} ({\Bbb R})$; in this case
$Hu=H_{0}u+ Vu$. It follows that, for $u\in{\cal D}(H)$, the function $u '''(x)-  v_{1}(x)u '(x)$ is absolutely continuous and
\[
(Hu)(x)=( u '''(x)-  v_{1}(x)u '(x))' + v_{0}(x)u  (x) .
\]

Let us discuss   main steps of a construction   of scattering theory for the operator $H$.    The resolvent  $R_{0}(z)= (H_{0}-z)^{-1}$, $z\in {\Bbb C}\setminus [0,\infty)$, of the operator $H_0 $ can be calculated   explicitly.

\begin{lemma}\label{2.1} 
Let $z\in{\Bbb C}\setminus [0,\infty)$ and $\zeta^4=z$, $\arg \zeta\in (0,\pi/2)$. Then
\begin{equation}
 (R_0(z)f)(x)=\frac{1}{4\zeta^3} \int_{-\infty}^\infty(  
 i e^{i\zeta |x-y|} -e^{-\zeta |x-y|} )f(y)dy.
\label{eq:R}\end{equation}
\end{lemma}

\begin{pf}
Using the Fourier   transform, we see that
\[
 (R_0(z)f)(x)=2\pi^{-1}\int_{-\infty}^\infty \int_{-\infty}^\infty  e^{ik (x-y)}(k^4-z)^{-1}  f(y)
dk   dy   .
\]
Therefore the resolvent kernel equals
\[
  R_0(x,y;z) =2\pi^{-1}\int_{-\infty}^\infty    e^{ik |x-y|}(k^4-z)^{-1}  dk  .
\]
This integral can be complemented in the upper   half-plane by a big half-circle and then calculated by residues at the points $\zeta$ and $i\zeta$.
\end{pf}

Let us denote by $\Pi$ the complex plane cut along $[0,\infty)$ (including upper and lower edges). According to Lemma~\ref{2.1} the resolvent kernel $ R_0(x,y;z)$ is a continuous function of $z\in \Pi$     with an exception of the point $z=0$.
The following two results are   immediate consequences of explicit  formula
 \e{eq:R}.
 
  \begin{proposition}\label{2.2n}
  If $f\in L_1 ({\Bbb R})$, then for all $z\in \Pi$, $z\neq 0$, the function $R_{0}(z)f\in C^3({\Bbb R})$ and $(R_{0}(z)f)'''(x)$ is an absolutely continuous function.
  \end{proposition}
 
 \begin{proposition}\label{2.2}
Let $G_j$ be the operator of multiplication by a function $g_j \in L_{2}({\Bbb R})$, $j=0, 1$.   Then
 the operator-valued functions
$$
 G_{0} D^l R_{0}(z) G_{1},   \q l=0,1, 2,3, 
 $$
 depend in the Hilbert-Schmidt norm continuously on $z\in \Pi$, $z\neq 0$.     
\end{proposition}

To extend the latter result to the resolvent $R (z) =(H-z)^{-1}$ of the operator $H$, we proceed from the   resolvent identity
\begin{equation}
R(z)=R_{0}(z)-R_{0}(z)VR (z)=R_{0}(z)-R (z)VR_{0} (z), \q \Im z\neq 0,
\label{eq:res}\end{equation}
where $V$ is operator \e{eq:V}.
Let $G_{j}$ and $\Omega_{j}$ be   operators of multiplication by the functions $\sqrt{ |v_{j}(x)}|$ and $\sgn v_{j}(x) $, respectively. We introduce an auxiliary space   ${\cal G}={\cal H}\oplus {\cal H}$ and define (bounded) operators ${\bf G}_{0}$, ${\bf G}: {\HH}^2 ({\Bbb R}) \ri {\cal G} $ by formulas\footnote{Although it is more convenient to write vectors as rows, we regard them as columns as far as matrix multiplication is concerned.}
\begin{equation}
{\bf G}_{0} =(\Omega_{0} G_{0}  , \Omega_{1} G_{1} D) , \q
{\bf G}  =(  G_{0}  ,   G_{1} D ).
\label{eq:G}\end{equation}
It follows from equality \e{eq:V} that $V= {\bf G}^* {\bf G}_{0}={\bf G}_{0}^* {\bf G}$. The resolvent identity
 \e{eq:res} implies that
 \[
(I+  {\bf G}_{0}R_{0}(z){\bf G}^* )(I-  {\bf G}_{0}R (z){\bf G}^* )=
(I-  {\bf G}_{0}R (z){\bf G}^* ) (I+  {\bf G}_{0}R_{0}(z){\bf G}^* )=I.
\]
Hence the inverse operator $(I+  {\bf G}_{0}R_{0}(z){\bf G}^* )^{-1}$ exists and is bounded
so that using again \e{eq:res} we obtain the representation
\begin{equation}
R(z)=R_{0}(z) -R_{0}(z) {\bf G}^* (I+  {\bf G}_{0}R_{0}(z){\bf G}^* )^{-1} {\bf G}_{0}R_{0} (z) .
\label{eq:res2}\end{equation}
Thus, the resolvent $R(z)$ for $\Im z\neq 0$   considered as  a mapping from ${\HH}^{-2} ({\Bbb R})$ to ${\HH}^2 ({\Bbb R})$ is a bounded operator.

It follows from Proposition~\ref{2.2} that under assumption \e{eq:2.1} the operator-valued function ${\bf G}_{0}R_{0}(z){\bf G}^*$, analytic for $\Im z\in {\Bbb C}\setminus [0,\infty)$, is continuous in the Hilbert-Schmidt norm  for $z\in \Pi$ except the point $z=0$. Therefore according to the analytic Fredholm alternative
(see, e.g., \cite{Ya}) the set ${\cal N}\in{\Bbb R}_{+}$  where at least one of two homogeneous equations
\begin{equation}
{\bf f}+  {\bf G}_{0}R_{0}(\lambda\pm i0){\bf G}^* {\bf f}=0, \q {\bf f}=(f_{0}, f_{1}),
\label{eq:hom}\end{equation}
has a non-trivial solution ${\bf f}\in {\cal G}$ is closed and has the Lebesgue measure zero. The operator-valued function $(I+  {\bf G}_{0}R_{0}(z){\bf G}^* )^{-1}$ of $z\in\Pi$ is continuous in norm   except   points from the set ${\cal N}\cup\{0\}$. Therefore equation \e{eq:res2} leads to the following result.

 \begin{theorem}\label{2.3} 
Let  assumption \e{eq:2.1} hold, and let $G_{j}$ be the same operators
as in Proposition~$\ref{2.2}$. Then
 the operator-valued functions $G_{1}  R (z) G_{2}$, $G_{1} D R (z) G_{2}$ and $G_{1}  D R (z) D G_{2}$  of $z\in\Pi$ are continuous in  the Hilbert-Schmidt norm   except   points from the set ${\cal N}\cup\{0\}$. The set $\Lambda={\Bbb R}_{+}\setminus {\cal N}$  is open and has   full Lebesgue measure.
\end{theorem}

 \begin{corollary}\label{2.3c} 
The spectrum of the operator $H$ on the set $\Lambda$ is absolutely continuous.
\end{corollary}

We denote by  $P^{(c)}$   the orthogonal projector on the absolutely continuous subspace ${\cal H}^{(c)}$ of the operator $H$.

\medskip

{\bf 2.2. }
Given Theorem~\ref{2.3}, an expansion in eigenfunctions of the operator $H$, a formula representation of the scattering matrix, etc.,  are consequences of general results of scattering theory (see \cite{BY4} or \cite{Ya}).

  Denote by $E_{0}(\lambda)$
the spectral family of the operator $H_{0}$. In the momentum representation the operator $E_{0}(\lambda)$ acts as multiplication by the characteristic function $\chi_\lambda(\xi)$ of the interval $(-\sqrt[4]{\lambda}, \sqrt[4]{\lambda})$, that is
\[
(\widehat{E_{0}(\lambda)f_{0}}) (\xi) =\chi_\lambda(\xi) \hat{f_{0}}(\xi)
\]
where $\hat{f_{0}}(\xi)$ is the Fourier transform of the function $f_{0}(x)$. It follows that for $f_{0}\in L_1  ({\Bbb R})\cap L_{2}  ({\Bbb R})$ the function $ (E_{0}(\lambda)f_{0},f_{0})$ belongs to the class 
$C^1 ({\Bbb R}_{+})$ and
\begin{equation}
d ( E_{0}(\lambda)f_{0}, f_{0}) / d\lambda
=4^{-1} \lambda^{-3/4} (|\hat{f_{0}}(\sqrt[4]{\lambda})|^2+ |\hat{f_{0}}(-\sqrt[4]{\lambda})|^2).
\label{eq:E1}\end{equation}
Next we construct the canonical spectral represention of the operator $H_{0}$.
We set
\begin{equation}
\Gamma_{0} (\lambda) f_{0}= 2^{-1}\lambda^{-3/8} (\hat{f_{0}}(\lambda^{1/4}), \hat{f_{0}}(-\lambda^{1/4}))  , \q \Gamma_{0} (\lambda): 
L_1  ({\Bbb R})\cap L_{2}  ({\Bbb R}) \ri {\Bbb C}^2,
\label{eq:spe}\end{equation} 
and define the  operator ${\cal F}_{0}: L_1  ({\Bbb R})\cap L_{2}  ({\Bbb R}) \ri L_{2}({\Bbb R}_{+}; {\Bbb C}^2)$  by the formula
$
({\cal F }_{0} f)(\lambda)\\ = \Gamma_{0} (\lambda) f .
$ 
This operator extends by continuity to a 
unitary mapping ${\cal F}_{0}:{\cal H}\ri L_{2}({\Bbb R}_{+}; {\Bbb C}^2)$.
Then ${\cal F }_{0}H_{0}=A{\cal F }_{0}$ where $A$ is the operator of multiplication by $\lambda$ in the space $ L_{2}({\Bbb R}_{+}; {\Bbb C}^2)$.

Now we discuss generalizations of these objects for the operator $H$.
Let us set
\begin{equation} 
  \Gamma_\pm(\lambda)f=\Gamma_0(\lambda)(I-VR(\lambda\pm i0))f,\quad f\in L_1  ({\Bbb R})\cap L_{2}  ({\Bbb R}),\q
\lambda\in \Lambda.
\label{eq:expsm}\end{equation}
According to Theorem~\ref{2.3}, $\Gamma_\pm(\lambda)f  $ is a   continuous function of $\lambda\in \Lambda$. 
A proof of the following result can be found in  \cite{BY4} or \cite{Ya}.

\begin{theorem}\label{EXPsm}
Let   assumption  \e{eq:2.1}   hold.   
Define the mapping ${\cal F}_\pm$ on the set $L_1  ({\Bbb R})\cap L_{2}  ({\Bbb R})$  by equalities $(\ref{eq:expsm})$ and 
\begin{equation}
({\cal F }_{\pm} f)(\lambda)= \Gamma_{\pm} (\lambda) f.  
\label{eq:expsm1}\end{equation}
 This mapping   extends by continuity to a
bounded operator ${\cal F}_{\pm}:{\cal H}\ri L_{2}({\Bbb R}_{+}; {\Bbb C}^2)$, satisfies the relations  
\begin{equation}
 {\cal F}_{\pm} {\cal F}_{\pm}^* =  I, \q {\cal F}_{\pm}^*{\cal F}_{\pm}=P^{(c)} 
 \label{eq:F1}\end{equation}
 and diagonalizes $H$, that is
 \begin{equation}
 {\cal F }_{\pm}H =A{\cal F }_\pm.
 \label{eq:F2}\end{equation}  
\end{theorem}
 
Time-dependent   wave operators for the pair $H_{0}$, $H$ are defined as strong limits
\begin{equation}
W_\pm=W_\pm(H,H_{0})  =\slim_{t\ri\pm\infty}e^{i H t} e^{-i H_{0} t}.
\label{eq:WO}\end{equation}
Recall that, by the spectral theorem,
\[
R (\lambda\pm i\varepsilon)=\pm \int_{0}^\infty e^{-\varepsilon t\pm i \lambda t} e^{\mp iHt} dt
\]
so that, by the Parseval identity, 
\begin{equation}
2\varepsilon\int_0^\infty e^{-2\varepsilon t}(e^{\mp itH_0}f_0,e^{\mp itH }f)dt=\pi^{-1}\varepsilon
\int_{-\infty}^\infty (R_0(\lambda\pm i\varepsilon)f_0, R (\lambda\pm i\varepsilon)f )d\lambda 
\label{eq:Parsid}\end{equation}
 for all $f_0,f\in{\cal H}$.
 
 Let us show that this expression   has a  limit as $\varepsilon\to 0$.
   This entails the existence of  the weak wave operators understood, moreover, in the Abelian sense. Such wave operators are defined by the limit of the left-hand side of \e{eq:Parsid}.  It suffices to verify the existence of the limit   for 
 $f_{0}, f \in L_1  ({\Bbb R})\cap L_{2}  ({\Bbb R})$.
 
 Let us consider the right-hand side of \e{eq:Parsid}. Note that, by the resolvent identity (\ref{eq:res}),
\[
\pi^{-1}\varepsilon  (R_0(\lambda\pm i\varepsilon)f_0, R (\lambda\pm i\varepsilon)f
)=(\delta_\varepsilon(H_0-\lambda)f_0, (I-V R (\lambda\pm i\varepsilon))f), 
\]
where
 \[
\delta_\varepsilon(H_0-\lambda)= (2\pi i)^{-1} (R_0(\lambda + i\varepsilon)
-R_0(\lambda - i\varepsilon))
\]
is an ``approximate" operator-valued delta-function. It follows from the spectral theorem, standard properties of the Cauchy type singular integrals and formulas \e{eq:E1}, \e{eq:spe} that
\begin{equation}
\lim_{\varepsilon\ri 0}(\delta_\varepsilon(H_0-\lambda) f_{0}, f )= d(E_{0}(\lambda)f_{0}, f ) / d\lambda =  \langle \Gamma_{0} (\lambda)f_0,   \Gamma_{0} (\lambda)  f \rangle 
\label{eq:Cauchy}\end{equation}
where ${ \langle} \cdot, \cdot { \rangle}$ is the scalar product in ${\Bbb C}^2$.
Therefore according to Proposition~\ref{2.2} we have 
\begin{equation}
\lim_{\varepsilon\rightarrow 0}{\bf G}_{0}\delta_\varepsilon(H_0-\lambda)f_{0}  
= {\bf G}_{0} \Gamma_0^\ast (\lambda)    \Gamma_0 (\lambda)f_{0} .
\label{eq:deltag}\end{equation}
The convergence in \e{eq:Cauchy} and \e{eq:deltag} is uniform on compact intervals of ${\Bbb R}_{+}$. Furthermore,   Theorem~\ref{2.3} ensures that there exists
\begin{equation}
\lim_{\varepsilon\rightarrow 0}{\bf G} R(\lambda\pm i\varepsilon) f. 
= {\bf G}  R(\lambda\pm i 0) f .
\label{eq:delt}\end{equation}
Combining relations \e{eq:Cauchy}, \e{eq:deltag} and \e{eq:delt}, we obtain  that
\begin{eqnarray}
\lim_{\varepsilon\ri 0}(\delta_\varepsilon(H_0-\lambda)f_0, (I-V R (\lambda\pm i\varepsilon))f)={ \langle}  \Gamma_{0} (\lambda)f_0,  \Gamma_0 (\lambda)f { \rangle} 
\nonumber\\
-(  {\bf G}_{0} \Gamma_0^\ast (\lambda)  \Gamma_{0} (\lambda)f_0,  {\bf G}
R (\lambda\pm i0) f )= { \langle}  \Gamma_{0} (\lambda)f_0,  \Gamma_{\pm} (\lambda)f
{ \rangle}.  
\label{eq:Parsid2}\end{eqnarray}
The convergence in \e{eq:delt} and \e{eq:Parsid2} is uniform with respect to $\lambda$ from compact intervals of the set $\Lambda$.

It remains to
 justify a passage to the limit $\varepsilon\ri 0$ in the integral in the right-hand side of (\ref{eq:Parsid}). Note that, by the Schwarz inequality, for any Borel set $X\subset {\Bbb R}$
 \begin{eqnarray}
   \Big|\int_X \varepsilon (R_0(\lambda\pm i\varepsilon)f_0, R (\lambda\pm i\varepsilon)f )d\lambda
\Big|^2
\nonumber\\
\leq\int_X \varepsilon \| R_0(\lambda\pm i\varepsilon)f_0\|^2d\lambda
\int_{-\infty}^\infty \varepsilon \| R (\lambda\pm i\varepsilon)f \|^2d\lambda
\nonumber\\
= \pi^2\int_X  (\delta_\varepsilon(H_0-\lambda)f_0,f_0)d\lambda\;\|f\|^2.
\label{eq:expwo4}\end{eqnarray}
Since the function $(\delta_\varepsilon(H_0-\lambda)f_0,f)$ is the Poisson integral of      function \e{eq:Cauchy} belonging to the space $L_{1}({\Bbb R})$, we see (see, e.g., \cite{Hof}) that the convergence in \e{eq:Cauchy} holds true
in the sense of $L_{1}({\Bbb R})$. Therefore
 the right-hand side of \e{eq:expwo4} tends to zero as $|X|\ri  0$ or
 $X=(N, \infty )$ and $N\ri \infty$ uniformly with respect to $\varepsilon\in (0,1)$.
 Moreover, it tends to zero  as $\varepsilon\ri 0$ if $X=(-\infty, 0)$.

 Thus, we have shown that the integral in the right-hand side of (\ref{eq:Parsid}) has the limit which equals the integral of function \e{eq:Parsid2} over $\lambda\in {\Bbb R}_{+}$. It follows     that there exists the limit
 \begin{equation}
\lim_{\varepsilon\to 0}2\varepsilon\int_0^\infty e^{-2\varepsilon t}(e^{\mp itH_0}f_0,e^{\mp itH }f)dt= 
\int_0^\infty  \langle  \Gamma_{0} (\lambda)f_0,  \Gamma_{\pm} (\lambda)f \rangle  d\lambda=({\cal F}_{0}f_0, {\cal F}_\pm f)
\label{eq:ParsidX}\end{equation}
and hence the Abelian weak wave operators for the pair $H_{0}$, $H$ exist and are equal to
${\cal F}_{\pm}^*  {\cal F}_{0}$. 

The strong wave operators \e{eq:WO}  also exist. This fact can be deduced from  general results of \cite{BY4} or \cite{Ya}. Alternatively, Theorem~\ref{2.3} entails that the operators ${\bf G}_{0}$ and ${\bf G}$ are $H$-smooth (as well as $H_{0}$-smooth)  in the sense of Kato (see, e.g., \cite{RS} or \cite{Ya}) on all compact intervals $X\subset \Lambda$ which also implies  the existence of  strong limits \e{eq:WO}. Finally, we note that under assumption \e{eq:2.1} the difference $R(z)-R_{0}(z) $ belongs to the trace class. Therefore the existence and completeness of wave operators \e{eq:WO} is a consequence of the   Birman-Kre\u{\i}n theorem obtained in \cite{BK}.

The results discussed above can be summarized in the following assertion.

 \begin{theorem}\label{2.4} 
 Let  assumption \e{eq:2.1} hold. 
Then wave operators 
\e{eq:WO} exist and satisfy the relation
$ W_{\pm} =  {\cal F}_{\pm}^*  {\cal F}_{0}$.   
The  operators $W_{\pm}$ are isometric and complete, that is their ranges $\ran W_\pm   ={\cal H}^{(c)}$. The intertwining property 
$ H W_\pm=W_\pm H_{0} $
holds. 
\end{theorem}

\medskip

{\bf 2.3. }
Since the scattering operator ${\cal S}=W_{+}^* W_{-}$ commutes  with the operator $H_{0}$, the operator 
\begin{equation}
{\cal F}_{0}{\cal S} {\cal F}_{0}^*= {\cal F}_{+}{\cal F}_{-}^*
\label{eq:S}\end{equation}
 acts in the space $L_{2}({\Bbb R}_{+}; {\Bbb C}^2)$ as multiplication by a $2\times 2$ matrix-valued function
\begin{equation}
S(\lambda)=\left(\begin{array}{cc}
s_{11}(\lambda)&s_{12}(\lambda)
\\
s_{21}(\lambda)&s_{22}(\lambda)
\end{array} \right) 
\label{eq:S1}\end{equation}
known as the scattering matrix. According to \e{eq:expsm1}
equality \e{eq:S} means that 
\begin{equation}
S(\lambda)\Gamma_{-}(\lambda)f =\Gamma_{+}(\lambda)f, \q \lambda\in \Lambda,
\q \forall f\in L_1({\Bbb R}) \cap  L_{2} ({\Bbb R}).
\label{eq:S5}\end{equation}

This equation determines the scattering matrix uniquely. Let us show that its solution is given by the formula
\begin{equation} 
 S(\lambda)= I-2\pi i \Gamma_{0}(\lambda) (V-VR(\lambda+i0) V)\Gamma_{0}^*(\lambda), \q \lambda\in \Lambda.
\label{eq:S3}\end{equation}
Let us set $\Gamma_{0} =\Gamma_{0}(\lambda)$, $R_{0}=R_{0}(\lambda+i0)$ and $R=R(\lambda+i0)$.
By virtue of definition \e{eq:expsm}, we have to check that
\[ 
 ( I-2\pi i \Gamma_{0}  (V-VR  V)\Gamma_{0}^*) \Gamma_{0}(I-VR^*)
 = \Gamma_{0}(I-VR). 
\]
In view of identity   \e{eq:Cauchy}, it suffices to verify  that
\[
(I-R  V)(R_{0}-R_{0}^*) (I-V R^*)
=R-R^*.
\]
This equality is a direct consequence of the resolvent identity \e{eq:res}.
Thus, we have proven the following result.

\begin{theorem}\label{SM}
Let  assumption  \e{eq:2.1}   hold.   
Then the scattering matrix $S(\lambda)$ for the pair $H_{0}$, $H$ admits representation \e{eq:S3}.  
\end{theorem}

According to \e{eq:S3}  the scattering matrix $S(\lambda)$ is a continuous function of $\lambda\in \Lambda$.

Since  $v_{j}(x)=\overline{v_{j}(x)}$ for $j=0,1$, the resolvent $R(z)$ commutes with the complex conjugation   which can  formally be written in terms of its kernel as
$\overline{R(x, y ;z) }=R(x, y ;\bar{z})$. 
Taking also into account that $R(\bar{z})=R^*(z)$, we see that the Green function $R(x,y ;z)$ is symmetric, that is
\begin{equation}
 R(x,y;z)=R(y ,x;z).
\label{eq:resc}\end{equation}
In view of representation \e{eq:S3}, it follows from this relation  that
\begin{equation} 
 s_{11}(\lambda)=s_{22}(\lambda).
\label{eq:S4}\end{equation}


\section{A homogeneous equation}
 
 
 {\bf 3.1. }
  Let us study a structure of the exceptional set ${\cal N}$.  Suppose   that, for one of the signs, equation  \e{eq:hom} is satisfied.  Set
\begin{equation}
\psi =R_{0}(\lambda \pm i0){\bf G}^* {\bf f}. 
\label{eq:hom1}\end{equation}
Then it follows from \e{eq:hom} that
\begin{equation}
{\bf f} + {\bf G}_{0}\psi=0 
\label{eq:hom2}\end{equation}
and hence
\begin{equation}
\psi + R_{0}(\lambda \pm i0)V \psi=0. 
\label{eq:hom3}\end{equation} 
    Taking into account definition  \e{eq:G} of the operator ${\bf G}$,  we  see that
 \begin{equation}
  V \psi    = -{\bf G}^*{\bf f} = \varphi_{0}+D \varphi_{1} 
\q
\mathrm{where } \q \varphi_{j}= -\sqrt{|v_{j}|} f_{j}   \in L_{1}({\Bbb R}),  \q   j=0,1.
\label{eq:V1}\end{equation}
The functions $R_{0}(\lambda \pm i0)\varphi_{0}$ and $R_{0}(\lambda \pm i0)D \varphi_1$ are well defined by Proposition~\ref{2.2n}.
  
  \begin{proposition}\label{SiDi} 
 Suppose that assumption  \e{eq:2.1}    holds.
 Let ${\bf f}\in {\cal G}$, ${\bf f}\neq 0$,   satisfy   equation  \e{eq:hom}, and let $\psi$ be defined by formula \e{eq:hom1}. Then $\psi$ is not identiacally zero, $\psi\in C^2 ({\Bbb R})$, the functions $\psi''(x)$ and 
 $\psi'''(x)-v_{1}(x) \psi'(x)$ are absolutely continuous and differential equation \e{eq:dex}
is satisfied.
\end{proposition}

\begin{pf}
 If $\psi = 0$, then ${\bf f}= 0$ according to equation \e{eq:hom2}.
 By virtue of Proposition~\ref{2.2n}, the inclusion  $\psi\in C^2 ({\Bbb R})$ and the absolute continuity of  $\psi''(x)$ 
follow from equation \e{eq:hom3} and representation \e{eq:V1}. For $\theta\in C_{0}^\infty ({\Bbb R})$, we have that
\[
(R_{0}(\lambda \pm i0)V \psi, (D^4-\lambda) \theta)=
 ( v_{1}\psi', R_{0}(\lambda \mp i0)(D^4-\lambda) \theta')+( v_{0}\psi, R_{0}(\lambda \mp i0)(D^4-\lambda) \theta) .
 \]
 Using equation \e{eq:hom3} and the fact that $R_{0}(\lambda \mp i0)(D^4-\lambda)\theta =\theta$ for all $\theta\in C_{0}^\infty ({\Bbb R})$, we can rewrite this equality as
 \[
- (  \psi, (D^4-\lambda) \theta) = ( v_{1}\psi',   \theta') + ( v_{0}\psi,   \theta).
 \]
 Thus, 
  \[
(-  \psi''' + v_{1}\psi',   \theta') =  (\lambda \psi- v_{0}\psi,   \theta)
 \]
 where $-  \psi''' + v_{1}\psi' \in L_{1}^{(loc)}({\Bbb R})$ and $\lambda \psi- v_{0}\psi\in L_{1}^{(loc)}({\Bbb R})$. It follows that the derivative in the sense of distributions of the function $  \psi''' - v_{1}\psi'$ equals $\lambda \psi- v_{0}\psi$. This implies that
 the function $  \psi''' - v_{1}\psi'$ is absolutely continuous and    differential equation \e{eq:dex} is satisfied. 
   \end{pf}
   
   Below differential equation \e{eq:dex} is always understood in the sense specified in Proposition~\ref{SiDi}. 
  
Next we  find asymptotic behavior of $\psi (x)$ as $|x|\to \infty$. To that end, we need   the following standard assertion.

\begin{lemma}\label{Si1} 
 Suppose   that, for one of the signs, equation  \e{eq:hom} is satisfied. Then
\begin{equation}
\Gamma_{0}(\lambda){\bf G}^* {\bf f} =0.
\label{eq:E2}\end{equation}
\end{lemma}

\begin{pf}
It follows from equation \e{eq:hom}   that
\[
\lim_{\varepsilon\ri 0}({\bf f}+  {\bf G}_{0}R_{0}(\lambda\pm i \varepsilon){\bf G}^* {\bf f}, {\bf G} R_{0}(\lambda\pm i \varepsilon){\bf G}^* {\bf f})=0.
\]
Taking here the imaginary part and using  that the operator $V={\bf G}^*{\bf G}_{0}$ is symmetric, we obtain the equality
\[
\lim_{\varepsilon\ri 0}( (R_{0}(\lambda+ i \varepsilon)-R_{0}(\lambda- i \varepsilon)){\bf G}^* {\bf f}, {\bf G}^* {\bf f})=0.
\] 
Now \e{eq:E2} is a consequence of   relation  \e{eq:Cauchy}.
\end{pf}

Below integrals containing derivatives of $L_{1}$-functions (for example, $ \varphi_{1}'(x)$) are understood in the sense of distributions, that is integration by parts is tacitly assumed. Using  formula
 \e{eq:E1} and representation \e{eq:V1}, we can reformulate Lemma~\ref{Si1} 
 in the following way.

\begin{corollary}\label{Si2} 
 Suppose   that, for one of the signs, equation  \e{eq:hom} is satisfied. Define the function $\psi$ by formula \e{eq:hom1}.
 Then for both signs $``\pm"$
\begin{equation}
\int_{-\infty}^\infty e^{\pm ikx} (V \psi )(x)  dx=0, \q k= \sqrt[4]{\lambda}.
\label{eq:E5}\end{equation}
\end{corollary}

We use also the following simple result.

\begin{lemma}\label{As} 
Let $\varphi \in L_{1}({\Bbb R})$. Then
\begin{equation}
( R_{0}(\lambda + i0) \varphi)(x) 
= \frac{i}{4 k^3}e^{\pm i k x} \int_{-\infty}^\infty e^{\mp i k  y}\varphi (y)dy
+ o(1) 
\label{eq:RR}\end{equation} 
as $x\ri \pm \infty$. 
\end{lemma}

\begin{pf}
 Suppose for definiteness
that  $x \ri + \infty$.  According to \e{eq:R} the function $ -  4 i k^3 ( R_{0}(\lambda + i0) \varphi)(x)$ consists of two terms. The first of them equals
\begin{eqnarray}
  \int_{-\infty}^\infty e^{i k |x-y|}\varphi(y)dy 
= e^{i k x} \int_{-\infty}^\infty e^{-i k  y}\varphi (y)dy
   \nonumber\\
 -e^{i k x} \int_x^\infty  e^{-i k  y} \varphi (y)dy
  + e^{- i k x} \int_x^\infty e^{i  ky}\varphi (y)dy.  
\label{eq:R1}\end{eqnarray}
Since $\varphi \in L_{1}({\Bbb R})$,  both integrals over $(x,\infty)$ in the right-hand side tend to zero as $x\ri + \infty$.   The second term equals
\begin{equation}
  \int_{-\infty}^\infty e^{- k |x-y|}\varphi(y)dy= e^{- k x} \int_{-\infty}^x e^{ k  y}\varphi(y)dy
  + e^{  k x} \int_x^\infty e^{- ky}\varphi(y)dy. 
\label{eq:R2}\end{equation}
  Since
\[
\big| \int_{-\infty}^x e^{ k  y}\varphi(y)dy\big| \leq e^{ k  x/2}
 \int_{-\infty}^{x/2}  | \varphi(y)| dy +
 e^{ k  x }
 \int_{x/2}^x  | \varphi(y)| dy
 \]
 and
  \[
\big| \int_x^\infty e^{- k  y}\varphi(y)dy\big|
 \leq e^{- k  x }  \int_x^\infty |\varphi(y)|dy,
 \] 
 both terms in the right-hand side of \e{eq:R2} tend to zero as $x\ri + \infty$.
 \end{pf}
 
Of course asymptotics of $ ( R_{0}(\lambda - i0) \varphi)(x) $
is obtained from \e{eq:RR} by the complex conjugation.

In view of Lemma~\ref{As} equation \e{eq:hom3} and condition  \e{eq:E5} imply   that
 \begin{equation}
 \lim_{|x |\ri   \infty}\psi(x)=0.
\label{eq:E5w}\end{equation} 
 Let us formulate the results obtained in the following intermediary assertion.
 
\begin{proposition}\label{inter} 
Under the  assumptions of Proposition~$\ref{SiDi}$ condition \e{eq:E5w} is satisfied.
\end{proposition}

 \medskip

 {\bf 3.2. }
In this subsection we use specific methods of ordinary differential equations.
  Let us first of all rewrite equation  \e{eq:dex} as a system of four equations of the first order. We set
  \begin{equation}
 {\bf u}=(u_{1} u_{2}, u_{3}, u_{4}), \q \mathrm{where}\q u_{1}=u,\; u_{2}=u',\;u_{3}=u'',\;
 u_{4}=u'''-v_{1}u',
\label{eq:syst}\end{equation}
and
\begin{equation} 
A= \left(\begin{array}{cccc}
 0&1&0&0
\\ 
0&0&1&0
\\ 
0&0&0&1
\\ 
\lambda&0&0&0
\end{array} \right)
,\q
K(x)= \left(\begin{array}{cccc}
 0&0&0&0
\\ 
0&0&0&0
\\ 
0&v_{1}(x)&0&0
\\ 
-v_{0}(x)&0&0&0
\end{array} \right).
\label{eq:syst1}\end{equation}
Then equation  \e{eq:dex} is equivalent to the system
 \begin{equation}
 {\bf u}'(x)=A  {\bf u}(x) + K(x) {\bf u}(x).
\label{eq:syst2}\end{equation}
Clearly, the matrix $A$ has eigenvalues $\gamma_{1}= ik$, $\gamma_{2}= -ik$,$\gamma_{3}= - k$,$\gamma_{4}= k$. We denote by ${\bf p}_{j}=(1, \gamma_{j},\gamma^2_{j},\gamma^3_{j})$, $j=1,2,3,4$, the corresponding eigenvectors. Let (non-orthogonal) projectors $P_{j}$ be defined by the relation
\[
P_{j} {\bf f}= c_{j}{\bf p}_{j} \q \mathrm{if} \q {\bf f}= \sum_{l=1}^4 c_{l}{\bf p}_{l}.
\]
Then $P_{j}^2=P_{j}$, $P_{j} P_{l}=0$ if $j\neq l$,  $ AP_{j}=\gamma_{j}P_{j}$ and
\[
I=\sum_{j=1}^4  P_{j}, \q e^{Ax}=\sum_{j=1}^4 e^{\gamma_{j}x} P_{j} .
\]

Although the following  result is a particular case of Problem~29, Chapter~3, of \cite{CL},
we give its proof for a completeness of our presentation.

 \begin{proposition}\label{syst} 
Let  assumption  \e{eq:2.1}   hold. Then, for each of the signs $``\pm"$,  system  \e{eq:syst2} 
has four solutions ${\bf u}_{j}^{(\pm)}  (x,\lambda )$, $j=1,2,3,4$, such that
\begin{equation}
{\bf u}_{j}^{(\pm)} (x , \lambda)=e^{\pm \gamma_{j}x} ( {\bf p}_{j}+o(1)) 
\label{eq:syst3}\end{equation}
as $x\ri \pm \infty$.  Estimates of the remainders in \e{eq:syst3} are uniform with respect to $\lambda$ from   compact subintervals of ${\Bbb R}_{+}$.
\end{proposition}

\begin{pf}
We suppose that $x\to +\infty$ and omit the upper index $``\pm"$.
Pick some $j=1,2,3,4$.
Let us set 
 \[
Y_{1}(x)=Y_{2}(x)=e^{-kx} P_{3}, \q Y_{3}(x)=0,
 \q Y_{4}(x)=\sum_{l=1}^3 e^{\gamma_{l}x} P_{l} 
\]
and
 \[
 Z_{j}(x)=e^{Ax} - Y_{j}(x).
\]
Remark that
 \begin{equation}
 {\pmb |}  Y _{j}(x) {\pmb |} \leq C e^{(\sigma_{j}-k)x}  \q\mathrm{for} \q x\geq 0  \q\mathrm{and} \q 
 {\pmb |} Z_{j}(x) {\pmb |} \leq C e^{\sigma_{j}x}  \q\mathrm{for} \q x\leq 0 
\label{eq:syst7}\end{equation}
where $   \sigma_{j}=\Re \gamma_{j}$. Let us choose a   number $a$ such that
\begin{equation}
 2 C \int_a ^\infty {\pmb |} K(y) {\pmb |} dy \leq 1
\label{eq:syst10}\end{equation}
and consider an integral equation
 \begin{equation}
{\bf u}_{j}  (x) =e^{\gamma_{j}x} {\bf p}_{j}+\int_{a}^x  Y_{j}(x-y) K(y) {\bf u}_{j}   (y) dy-
\int_x ^\infty Z_{j}(x-y) K(y) {\bf u}_{j}   (y) dy.
\label{eq:syst6}\end{equation}

Below we also omit the index $j$. Let us show that equation \e{eq:syst6} has a solution ${\bf u} (x)$ satisfying an estimate
 \begin{equation}
{\pmb |} {\bf u} (x) {\pmb |} \leq 2 p e^{\sigma x},\q p={\bf p}, \q x\geq a.
\label{eq:syst8}\end{equation}
We use the method of successive approximations setting ${\bf u}^{(0)} (x)=  {\bf p} e^{\gamma x}$ and 
 \begin{equation}
{\bf u}^{(l+1)} (x) =  {\bf p} e^{\gamma x} +\int_{a}^x  Y (x-y) K(y) {\bf u}^{(l )} (y) dy-
\int_x ^\infty Z (x-y) K(y)  {\bf u}^{(l )} (y) dy.
\label{eq:syst9}\end{equation}
Let us check that, for all $l$,
 \begin{equation}
{\pmb |} {\bf u}^{(l )} (x) -{\bf u}^{(l -1)} (x) {\pmb |}
\leq 2^{-l } p e^{\sigma x}.
\label{eq:syst11}\end{equation}
Supposing \e{eq:syst11} for some $l$ and using definition \e{eq:syst9}, we obtain an estimate 
\[
{\pmb |} {\bf u}^{(l+1)} (x) -{\bf u}^{(l  )} (x){\pmb |} \leq 2^{-l } p \big(\int_{a}^x  {\pmb |} Y (x-y) {\pmb |}   {\pmb |}  K(y) {\pmb |}  e^{\sigma y} dy +
\int_x ^\infty{\pmb |} Z (x-y) {\pmb |}   {\pmb |}  K(y) {\pmb |}
e^{\sigma y} dy \big).
\]
By virtue of inequalities \e{eq:syst7} and condition \e{eq:syst10} this expression does not exceed
\[
2^{-l } p \: C e^{\sigma x} \int_{a}^\infty   {\pmb |}  K(y) {\pmb |} dy \leq 2^{-l -1} p e^{\sigma x}  .
\]
This proves estimate \e{eq:syst11} for $l+1$ in place of $l$ and hence for all $l$. Thus,   the sequence ${\bf u}^{(l )} (x) $ converges as $l\to\infty$ to a function  ${\bf u} (x) $ satisfying   estimate \e{eq:syst8}.  
Passing in \e{eq:syst9} to the limit $l\to\infty$, we get equation \e{eq:syst6}.

To prove asymptotics \e{eq:syst3} for the  function  ${\bf u} (x) $, we combine inequalities \e{eq:syst7} and \e{eq:syst8}. Obviously, the last integral in the right-hand side of \e{eq:syst6} is $o(e^{\sigma x})$ as $x\to\infty$. The first integral is estimated by
\[
2 C \: p e^{\sigma x} \big(  e^{- k x/2}\int_{a}^{x/2} {\pmb |}  K(y){\pmb |}dy  +
 \int_{x/2}^x {\pmb |} K(y) {\pmb |} dy   \big) 
 \]
 which is also $o(e^{\sigma x})$.
 
 Finally, a direct differentiation shows that a solution of integral equation \e{eq:syst6} 
satisfies also system \e{eq:syst2}.
 \end{pf}

 Let  $u _{j}^{(\pm)} (x , \lambda) $ be the first component of the vector ${\bf u} _{j}^{(\pm)} (x , \lambda) $. Proposition~\ref{syst} can be reformulated in terms of solutions of equation \e{eq:dex}.

 \begin{proposition}\label{CL} 
Let  assumption  \e{eq:2.1}   hold. Then, for each of the signs $``\pm"$,    differential equation \e{eq:dex} 
has four solutions $u_{j}^{(\pm)}(x, \lambda)$, $j=1,2,3,4$, such that
\begin{eqnarray}
u_{1}^{(\pm)}(x , \lambda)=e^{\pm ikx} (1+o(1)) ,\q u_{2}^{(\pm)}(x ,\lambda)=e^{\mp ikx}(1+o(1)),
\nonumber\\
 u_{3}^{(\pm)}(x , \lambda)=e^{\mp kx} (1+o(1)),\q u_{4}^{(\pm)}(x , \lambda)=e^{\pm kx}(1+o(1))
\label{eq:de1}\end{eqnarray}
as $x\ri\pm \infty$.  Estimates of the remainders in \e{eq:de1} are uniform with respect to $\lambda$ from   compact subintervals of ${\Bbb R}_{+}$.
\end{proposition}

 \begin{remark}\label{diff} 
 It follows from formulas  \e{eq:syst} that asymptotic relations \e{eq:de1} are two  times differentiable with respect to $x$. Moreover, 
 \[
 d^3 u_{j}^{(\pm)}(x , \lambda)/dx^3- v_{1}(x)d  u_{j}^{(\pm)}(x , \lambda)/dx
 =\pm \gamma_{j}^3 e^{\pm\gamma_{j}x} ( {\bf p}_{j}+o(1))
 \]
 as $x\to\pm\infty$. 
  \end{remark}

Every solution $u(x)$ of equation \e{eq:dex} is a linear combination of the solutions
  $u_{j}^{(\pm)}(x , \lambda)$, $j=1,2,3,4$. Therefore if $u(x)\to 0$   as $x\ri\pm\infty$, then    $u(x)= c^{(\pm)} u_{3}^{(\pm)}(x , \lambda)$
     for some constant $c^{(\pm)}$
and hence belongs to $L_{2}({\Bbb R}_{\pm})$. In particular, we have 

\begin{proposition}\label{CL1} 
Suppose that a function $u(x)$ satisfies   equation
\e{eq:dex} and $u (x) =o(1)$ as $x \ri \pm \infty$. Then $u (x) =0$ if $\lambda$ is not an eigenvalue  of the operator $H$. 
\end{proposition}

Combining Propositions~\ref{SiDi}, \ref{inter} and \ref{CL1}, we obtain  that every
$\lambda\in {\cal N}$ is necessarily an eigenvalue of the operator $H$. Taking also into account Theorem~\ref{2.4}, we see that the singular continuous spectrum of the operator $H$ is empty.
Conversely, if $\psi$ is an eigenfunction of the operator $H$, then it satisfies equation \e{eq:hom3} and hence ${\bf f}$ defined by \e{eq:hom2} satisfies equation \e{eq:hom}.
 Positive eigenvalues of   $H$  are of course simple.

Let us finally show that eigenvalues of the operator $H$ might accumulate at the point zero only. Suppose on the contrary that eigenvalues $\lambda_{n}=k_{n}^4\ri \lambda_{0}=k_{0}^4 >0$. Let $\psi_{n}$ be the corresponding normalized eigenfunctions.
By Proposition~\ref{CL}, we have that
\[
\psi_{n}(x)=a^{(\pm)}_{n}e^{-k_{n}| x |}(1+o(1))
\]
as $x\ri \pm \infty$. The estimate of the remainder here is uniform with respect to $n$. Since $\| \psi_{n}\|=1$, we have that $|a^{(\pm)}_{n}|\leq C<\infty$.  Therefore, for all $\varepsilon> 0$, and sufficiently large $R=R(\varepsilon)$
\begin{equation}
\int_{|x|\geq R}|\psi_{n}(x)|^2 dx<\varepsilon
\label{eq:comp}\end{equation}
 uniformly with respect to $n$. Moreover, we have that $h[\psi_{n}, \psi_{n}]=\lambda_{n}$ and  hence   $\| \psi_{n}\|_{{\HH }^2({\Bbb R})}\leq C<\infty$ according to estimate \e{eq:QF}. Together with \e{eq:comp}, this ensures compactness of the set of the functions $\psi_{n}$ in the space $L_{2}({\Bbb R})$ which contradicts their orthogonality.

Thus, we have obtained

\begin{theorem}\label{sing} 
Let  assumption \e{eq:2.1}  hold. Then  ${\cal N}$ coincides with the set  of positive eigenvalues of the operator $H$, the singular continuous spectrum of the operator $H$ is empty  and eigenvalues of the operator $H$ might accumulate at the point zero only.
\end{theorem}

According to Theorem~\ref{sing}  the absolutely continuous subspace ${\cal H}^{(c)}$ of the operator $H$ equals ${\cal H}^{(c)}={\cal H}\ominus {\cal H}^{(p)}$ where
 ${\cal H}^{(p)}$ is the subspace spanned by eigenfunctions (of the point spectrum) of the operator $H$.


\section{Eigenfunctions of the continuous spectrum}
 

{\bf 4.1. }
Eigenfunctions $\psi^{(\pm)}_j(x , \lambda)$, $j=1,2$,  of the continuous spectrum of the operator $H$ are defined by the formula
\begin{equation}
\psi^{(\pm)}_{j}(\lambda)= \psi^{(0)}_{j}(\lambda) -R(\lambda\mp i 0) V \psi^{(0)}_{j}(\lambda) , \q j=1,2,
\label{eq:WF}\end{equation}
where $\psi^{(0)}_{1}(x, \lambda)=e^{ikx}$, $\psi^{(0)}_2(x, \lambda)=e^{-ikx}$ and $k=\sqrt[4]{\lambda} > 0$. 
Set $\psi _{j}(x, \lambda)=\psi^{(-)}_{j}(x, \lambda)$, $j=1,2$. By virtue of  property
  \e{eq:resc} we have that $\psi^{(+)}_1(x, \lambda)=\overline{\psi _2(x, \lambda)}$ and $\psi^{(+)}_2(x, \lambda)=\overline{\psi _1(x,\lambda)}$.  
   The functions $\psi_{j}(x,\lambda)$ are known also as wave functions.
    The resolvent identity \e{eq:res} implies that the functions
$\psi _{j}^{(\pm)}(x, \lambda)$ satisfy the Lippmann-Schwinger equation
\begin{equation}
\psi_{j}^{(\pm)}(\lambda)= \psi^{(0)}_{j}(\lambda) -R_{0}(\lambda \mp  i 0) V \psi _{j}^{(\pm)}(\lambda) , \q j=1,2.
\label{eq:WF2}\end{equation}
Remark that the right-hand side here is correctly defined.  Indeed, it follows from Theorem~\ref{2.3} that  $g \psi _{j}(  \lambda), g \psi _{j}'(  \lambda) \in L_{2}({\Bbb R})$   for an arbitrary function $g\in L_{2}({\Bbb R})$ and hence the functions $ V \psi _{j}(  \lambda) $ admit representation \e{eq:V1}.

Similarly to the proof of of Proposition~\ref{SiDi}, it is easy to deduce from \e{eq:WF2} that the wave functions satisfy also differential equation \e{eq:dex}. Their asymptotics 
as $|x|\ri\infty$ can be found with a help of Lemma~\ref{As}.   Thus, we obtain the following result.

\begin{proposition}\label{WFS}
Let   assumption   \e{eq:2.1}  hold.   
Suppose that  $\lambda=k^4$ is not an eigenvalue of the operator $H$.
Let   matrix  \e{eq:S1} be defined by equation
\e{eq:S3}, and let  the  solutions 
$\psi_1(x,\lambda)$ and $\psi_2(x,\lambda)$ of equation $(\ref{eq:dex})$ be defined by  formula  $(\ref{eq:WF})$. Then the asymptotic relations
\begin{equation}\left\{\begin{array}{lcl} 
   \psi_1(x,\lambda) =e^{ik x}+s_{21}(\lambda) e^{-ik x}+o(1) ,\quad x\ri -\infty,
   \\
\psi_1(x,\lambda) = s_{11}(\lambda) e^{ik x}+o(1) ,\quad x\ri \infty,
\end{array}\right.
\label{eq:WFS1}\end{equation}
and
\begin{equation}\left\{\begin{array}{lcl} 
\psi_2(x,\lambda) =e^{-ik x}+s_{12}(\lambda) e^{ik x}+o(1) ,\quad x\ri \infty,
\\  
\psi_2(x,\lambda) = s_{22}(\lambda) e^{-ik x}+o(1) ,\quad x\ri -\infty,
\end{array}\right.
\label{eq:WFS2}\end{equation} 
hold.
\end{proposition}

\begin{remark}\label{di1}
Here and below 	all asymptotic relations are differentiable in the sense of Remark~\ref{diff}. Actually, for example, the first relation  \e{eq:WFS1} entails, by virtue of Proposition~\ref{CL}, that
\[
  \psi_1(x,\lambda)= s_{21}(\lambda) u_1^{(-)}(x,\lambda)+u_2^{(-)}(x,\lambda),
  \]
  where the functions $u_1^{(-)}$ and    $u_2^{(-)}$ are differentiable
according to  Remark~\ref{diff}.
  \end{remark}

Theorem~\ref{EXPsm} can be reformulated as an expansion of an arbitrary function in a generalized Fourier integral over the eigenfunctions $\psi_1^{(\pm)}(x,\lambda)$ and $\psi_2^{(\pm)}(x,\lambda)$ of the operator $H$. Indeed, according to definition   \e{eq:WF} for an arbitrary $f\in L_{1}({\Bbb R})\cap L_2({\Bbb R})$
\[
((I-V R(\lambda\pm i0))f, \psi_{j}^{(0)}(\lambda))= ( f,  \psi_{j}^{(\pm)}(\lambda)), \q j=1,2,
\]
and according to definitions  \e{eq:spe} and \e{eq:expsm}
\begin{equation}
\Gamma_{\pm} (\lambda) f= 2^{-1}(2\pi)^{-1/2} \lambda^{-3/8}
( ( f,  \psi_{1}^{(\pm)}(\lambda)) ,( f,  \psi_{2}^{(\pm)}(\lambda))).
\label{eq:FF1}\end{equation} 
Therefore  formula \e{eq:expsm1}  reads as $({\cal F }_{\pm} f)(\lambda) =(\tilde{f}_1(\lambda), \tilde{f}_2(\lambda))$ where
\[
\tilde{f}_j(\lambda)= 2^{-1}(2\pi)^{-1/2} \lambda^{-3/8}
 \int_{-\infty}^\infty \overline{\psi_{j}^{(\pm)}(x,\lambda)}f(x) dx . 
\]
It follows that the  relation  $P^{(c)}f={\cal F }_{\pm}^* {\cal F }_{\pm}f$ can (formally) be written as 
\[
(P^{(c)} f)(x)= 2^{-1}(2\pi)^{-1/2}   \sum_{j=1}^2 \int_0^\infty  \psi_j^{(\pm)}(x,\lambda)\tilde{f}_{j}(\lambda)
\lambda^{-3/8} d\lambda .
\]

By virtue of representation \e{eq:FF1},   equality \e{eq:S5}  is equivalent to relations
\begin{equation}\left\{\begin{array}{lcl}
 s_{11}(\lambda) \overline{\psi_1(x, \lambda)}+s_{12}(\lambda) \overline{\psi_2(x, \lambda)}&=&\psi_2(x, \lambda),
\\
s_{21}(\lambda) \overline{\psi_1(x, \lambda)}+s_{22}(\lambda) \overline{\psi_2(x, \lambda)}&=&\psi_1(x, \lambda).
\end{array}\right.
\label{eq:S2}\end{equation}

\medskip

{\bf 4.2.}
Let us show that asymptotics \e{eq:WFS1} or \e{eq:WFS2} determine uniquely solutions of equation (\ref{eq:dex}). We start with an auxiliary assertion which is true without assumption   \e{eq:2.1}.

\begin{lemma}\label{4di} 
Suppose that the functions $v_{0}$ and $v_{1}$ are real. Set
\begin{equation}
F_{u}(r)=     
(u^{\prime\prime\prime}(r)-v_{1}(r) u ^{\prime  }(r)) \overline{u (r)}
- u ^{\prime\prime }(r)\overline{ u ^{\prime}(r)} .
\label{eq:IP}\end{equation}
Then,
for an arbitrary solution $u(x)$ of equation \e{eq:dex} and all $r$, we have
\begin{equation}
\Im F_{u}(r)=\Im F_{u}(-r).
\label{eq:IP1}\end{equation}
\end{lemma}

\begin{pf}
It follows from equation (\ref{eq:dex})   that
\[
\Im \int_{-r}^r (u ''' -  v_{1}u^{\prime} )^{\prime} \overline{u}dx=0.
\]
Integrating here by parts, we get equality \e{eq:IP1}.
 \end{pf}

Now we can formulate the uniqueness result.

\begin{proposition}\label{WFSU}
Let   assumption   \e{eq:2.1}  hold.   
Suppose that  $\lambda=k^4$ is not an eigenvalue of the operator $H$.
If a  solution  
$u (x )$   of equation $(\ref{eq:dex})$ satisfies  the conditions
\begin{equation}\left\{\begin{array}{lcl} 
u (x ) = \sigma_{+} e^{\pm ik x}+o(1) ,\quad x\ri \infty,
\\  
u (x ) = \sigma_{-} e^{\mp ik x}+o(1) ,\quad x\ri -\infty,
\end{array}\right.
\label{eq:WFS3}\end{equation}
for one of the signs 
and  some numbers $\sigma_+$ and $\sigma_{-}$, then  $\sigma_+=\sigma_{-}=0$ and $u (x ) =0$.
\end{proposition}

\begin{pf}
It follows (see Remark~\ref{di1}) from  relations \e{eq:WFS3} that function \e{eq:IP} has asymptotics
\[
F_{u}(r)=  2 (\pm ik)^3 |\sigma_{+}|^2+ o(1), \q r\ri\infty,
\]
and
\[
F_{u}(-r)=   2 (\mp ik)^3 |\sigma_{-}|^2+ o(1), \q r\ri -\infty.
\]
Using Lemma~\ref{4di}, we find that $|\sigma_{+}|^2+ |\sigma_{-}|^2=0$ and hence $u (x)\ri 0$ as $|x|\ri \infty$. Thus,      $u (x)= 0$ by Proposition~\ref{CL1}.
 \end{pf}
 
 \begin{corollary}\label{WFSU1} 
If a  solution  
$\tilde{\psi}_{1}(x )$   of equation $(\ref{eq:dex})$ has asymptotics
\e{eq:WFS1} with some coefficients $s_{11}$ and $s_{21}$, then necessarily $s_{11}$ and $s_{21}$ are the entries of the scattering matrix and $\tilde{\psi}_{1}=\psi_{1}$. Similarly, if a  solution  
$\tilde{\psi}_{2}(x )$   of equation $(\ref{eq:dex})$ has asymptotics
\e{eq:WFS2} with some coefficients $s_{12}$ and $s_{22}$, then necessarily $s_{12}$ and $s_{22}$ are the entries of the scattering matrix and $\tilde{\psi}_{2}=\psi_{2}$.
\end{corollary}

 Formulas \e{eq:WFS1} and \e{eq:WFS2} give us the definition of the scattering matrix in terms of solutions of differential equation \e{eq:dex}. Similarly to the Schr\"odinger equation, the numbers $s_{11}(\lambda)$, $s_{22}(\lambda)$ can be
interpreted as transmission coefficients and  $s_{21}(\lambda)$, $s_{12}(\lambda)$ can be interpreted as reflection coefficients for a plane wave coming from minus or plus infinity and interacting with the potentials $v_{0}(x)$ and $v_{1}(x)$.

\medskip

{\bf 4.3.}
 Here we find asymptotics of the wave functions up to terms decaying super-exponentially at infinity. We say that a function $\varphi(x)$ decays super-exponentially if
\begin{equation}
\varphi(x)=O (e^{-\varkappa |x|}), \q \forall \varkappa>0, \q |x |\ri\infty.
\label{eq:super}\end{equation}
The following result supplements Lemma~\ref{As}.

\begin{lemma}\label{Asc} 
Let a function $\varphi  $ satisfy condition \e{eq:super}. Then
\begin{eqnarray*}
( R_{0}(\lambda + i0) \varphi)(x) 
&=& \frac{i}{4 k^3}e^{\pm i k x} \int_{-\infty}^\infty e^{\mp i k  y}\varphi (y)dy
\nonumber\\
&-&  \frac{1}{4 k^3}e^{- k | x|} \int_{-\infty}^\infty e^{\pm  k  y}\varphi (y)dy
+ O (e^{-\varkappa |x|}) , \q \forall \varkappa>0,
\end{eqnarray*}  
as $x\ri \pm \infty$.
\end{lemma}

\begin{pf}
 Suppose again for definiteness
that  $x \ri + \infty$.  According to \e{eq:R}  the function $ -4 i k^3 ( R_{0}(\lambda + i0) \varphi)(x)$ consists of the   terms \e{eq:R1} and 
\begin{eqnarray}
  \int_{-\infty}^\infty e^{- k |x-y|}\varphi(y)dy&=& e^{- k x} \int_{-\infty}^\infty e^{ k  y}\varphi(y)dy 
  \nonumber\\
  &-& e^{- k x} \int_{x}^\infty e^{ k  y}\varphi(y)dy
  + e^{  k x} \int_x^\infty e^{- ky}\varphi(y)dy.  
\label{eq:R2E}\end{eqnarray}
The integrals over $(x,\infty)$ in the right-hand sides of \e{eq:R1} and \e{eq:R2E}   decay   super-exponentially by virtue of condition \e{eq:super}. 
 \end{pf}
 
 Let us return to the Lippmann-Schwinger equation \e{eq:WF2} and take into account that the functions $\psi_j(x,\lambda)$ and $\psi_j'(x,\lambda)$, $j=1,2$, are bounded.
 Lemma~\ref{Asc} yields now a more precise form of Proposition~\ref{WFS}.
 
 \begin{proposition}\label{WFSE}
Let the functions $v_{0}$  and $v_{1}$ satisfy condition \e{eq:super}.
Suppose that  $\lambda=k^4$ is not an eigenvalue of the operator $H$.
Let   matrix  \e{eq:S1} be defined by equation
\e{eq:S3}, and let  the  solutions 
$\psi_1(x,\lambda)$ and $\psi_2(x,\lambda)$ of equation $(\ref{eq:dex})$ be defined by  formula  $(\ref{eq:WF})$. Then we have the asymptotic relations
\begin{equation}\left\{\begin{array}{lcl}
\psi_1(x,\lambda) = s_{11}(\lambda) e^{ik x}
+b_{11}(\lambda) e^{-k x}
   +O(e^{-\varkappa x})   ,\quad x\ri \infty,
   \\ 
   \psi_1(x,\lambda) =e^{ik x}+s_{21}(\lambda) e^{-ik x}+b_{21}(\lambda) e^{k x}
   +O(e^{\varkappa x}) ,\quad x\ri -\infty,
   \end{array}\right.
\label{eq:WFS1E}\end{equation}
and
\begin{equation}\left\{\begin{array}{lcl} 
\psi_2(x,\lambda) =e^{-ik x}+s_{12}(\lambda) e^{ik x} 
+b_{12}(\lambda) e^{-k x}
   +O(e^{-\varkappa x})  ,\quad x\ri \infty,
\\  
\psi_2(x,\lambda) = s_{22}(\lambda) e^{-ik x} +b_{22}(\lambda) e^{k x}
   +O(e^{\varkappa x})  ,\quad x\ri -\infty,
\end{array}\right. 
\label{eq:WFS2E}\end{equation} 
where  $\varkappa$ is arbitrary and
\[
b_{1l}(\lambda)=\frac{1}{4k^3}\int_{-\infty}^\infty e^{ky} (V  \psi_{l}(\lambda))(y)dy,
\q
b_{2l}(\lambda)=\frac{1}{4k^3}\int_{-\infty}^\infty e^{-ky} (V \psi_{l}(\lambda))(y)dy,
\]
$l=1,2$.
The asymptotic coefficients $b_{jl}(\lambda)$, $j, l =1,2$,  are continuous functions of $\lambda\in \Lambda$.
\end{proposition}

Let us find a relation between the sets of coefficients $s_{jl}(\lambda)$ and $b_{jl}(\lambda)$. To that end, we plug asymptotics \e{eq:WFS1E} and \e{eq:WFS2E}     into system \e{eq:S2}. Neglecting super-exponentially decaying terms,
 we have that, as $x\ri\infty$,
  \begin{equation} 
   s_{11} ( \overline{s_{11}}e^{-ikx}+  \overline{b_{11}}e^{-kx})
   + s_{12}  (   \overline{s_{12}}e^{-ikx}+  \overline{b_{12}}e^{-kx})
   = e^{-ik x}+ b_{12}  e^{-k x}
\label{eq:fourth}\end{equation}
   and
     \begin{eqnarray} 
   s_{21} ( \overline{s_{11}}e^{-ikx}+  \overline{b_{11}}e^{-kx})
   +s_{22} (e^{ikx}+  \overline{s_{12}}e^{-ikx}+  \overline{b_{12}}e^{-kx})
  \nonumber\\
   =  s_{11}e^{ikx}+   b_{11}e^{-kx}.
\label{eq:fourth1}\end{eqnarray}
 Similarly,   if $x\ri-\infty$, we have that
     \begin{eqnarray} 
  s_{11} (e^{-ikx} + \overline{s_{21}}e^{ikx}+  \overline{b_{21}}e^{kx})
   +s_{12} (   \overline{s_{22}}e^{ ikx}+  \overline{b_{22}}e^{ kx})
     \nonumber\\
   =  s_{22}e^{-ikx}+   b_{22}e^{kx}
\label{eq:fourth3}\end{eqnarray}
 and  
    \begin{equation}
   s_{21} (\overline{s_{21}}e^{ikx}+  \overline{b_{21}}e^{kx})
   +s_{22} (   \overline{s_{22}}e^{ ikx}+  \overline{b_{22}}e^{kx})
   = e^{ik x} +b_{21}  e^{k x}.
\label{eq:fourth2}\end{equation}

 Comparing the coefficients at $e^{-ik x}$ in the left- and right-hand sides of equations (\ref{eq:fourth}) and \e{eq:fourth1},
 we find that
\begin{equation}
 |s_{11}|^2 +|s_{12}|^2 =1
 \label{eq:U1}\end{equation}
 and
 \begin{equation}
 s_{21} \overline{s_{11}} + s_{22} \overline{s_{12}}=0.
 \label{eq:U2}\end{equation}
  Comparing the coefficients at $e^{ik x}$ in the left- and right-hand sides of equations (\ref{eq:fourth3}) and \e{eq:fourth2},
 we obtain again relations \e{eq:U2} and
\begin{equation}
 |s_{21}|^2 +|s_{22}|^2 =1.
 \label{eq:U3}\end{equation}
 Identities \e{eq:U1}, \e{eq:U2} and \e{eq:U3} show that the scattering matrix is a unitary operator in ${\Bbb C}^2$. This result has   already been obtained in Section~2 as a consequence of the completeness of the wave operators. 
 
 Comparing the coefficients at $e^{ik x}$ (or at $e^{-ik x}$) in the left- and right-hand sides of equation (\ref{eq:fourth1}) (or of equation (\ref{eq:fourth3})),
 we recover relation  \e{eq:S4}, which has already  been  obtained in Section~2 as a consequence of  the stationary representation \e{eq:S3} and of relation \e{eq:resc}.

Finally, comparing the coefficients at exponentially decreasing terms in equations 
 \e{eq:fourth} - \e{eq:fourth2}, we obtain four relations
between the coefficients $s_{jl} $ and $b_{jl} $. We write them in the matrix form as
\begin{equation}  \left(\begin{array}{cc}
 s_{11}(\lambda)&s_{12}(\lambda)
\\ s_{21}(\lambda)&s_{22}(\lambda)
\end{array} \right)\,
  \left(\begin{array}{cc} \overline{b_{11}}(\lambda)&
\overline{b_{21}}(\lambda)
\\ 
\overline{b_{12}}(\lambda)&
\overline{b_{22}}(\lambda)
\end{array} \right)= 
 \left(\begin{array}{cc} b_{12}(\lambda)&b_{22}(\lambda)
\\ b_{11}(\lambda)&b_{21}(\lambda)
\end{array} \right).
\label{eq:4di7}\end{equation}
We formulate this result in a following assertion.

\begin{theorem}\label{4di3}
Let the functions $v_{0}$  and $v_{1}$ satisfy condition \e{eq:super}.
Suppose that  $\lambda=k^4$ is not an eigenvalue of the operator $H$.
Then the asymptotic coefficients $s_{jl} $ and $b_{jl} $ in \e{eq:WFS1E} and
  \e{eq:WFS2E} are linked by relation  \e{eq:4di7}. 
\end{theorem}

We emphasize that relations  \e{eq:4di7} as well as \e{eq:S4} are consequences of the invariance of the problem with respect to the complex conjugation.

 \medskip 

{\bf 4.4.} 
The approach developed  for the problem on the whole line works of course  for a similar problem in the space $L_{2}({\Bbb R}_{+})$. In this case we have to add a boundary condition at the point $x=0$ which we choose as
\begin{equation}
u(0)=u^{\prime }(0)=0.
\label{eq:4ordfr}\end{equation}
The operator $H$ is defined by differential expression \e{eq:2.1a} on functions satisfying
\e{eq:4ordfr}. As a ``free" operator,  we
take    $H_0=D^4$ with the same boundary condition.
  Since the operator $H$ has a simple spectrum, there is now only one wave function $\psi(x , \lambda)$ and the scattering operator  acts in the space $L_{2}({\Bbb R}_{+})$ as multiplication by  the function (scattering matrix) $s(\lambda)$.

Let us construct the   function $\psi(x, \lambda)$. Remark first that the function
 \begin{equation}
  \psi_{0}(x,\lambda)= \cos ( kx+ \pi/4) - 2^{-1/2}  e^{-k x} 
\label{eq:5ordfr}\end{equation}
satisfies the equation $u^{(4)}(x)=\lambda u(x)$ and boundary condition \e{eq:4ordfr}.
Assume as usual  that
$\lambda$ is not an eigenvalue of $H$.  A solution $\psi(x, \lambda)$ of equation \e{eq:dex} is then defined (cf. \e{eq:WF})  by the formula
 $\psi( \lambda)=\psi_{0}( \lambda)-R(\lambda+ i0) V \psi_{0}( \lambda)$. Similarly to Proposition~\ref{WFS}, it can be shown under short-range assumption \e{eq:2.1} that
 \begin{equation}
\psi(x , \lambda)=  2^{-1}(s(\lambda)e^{ i kx+ \pi i/4}+ e^{-i kx-\pi i/4})+ o(1)
\label{eq:HL}\end{equation}
as $x\to\infty$.  Then we use (cf. Lemma~\ref{4di}) that $\Im F_u(r)=0$ for all solutions $u(x)$ of equation \e{eq:dex} satisfying \e{eq:4ordfr}. Applying this result to the difference of two functions obeying \e{eq:HL}, we see that actually condition \e{eq:HL} distinguishes a unique solution of \e{eq:dex}. In particular,
  \begin{equation}
\overline{\psi(x , \lambda) }=   s(\lambda)\psi(x , \lambda).
\label{eq:uni}\end{equation} 
  Moreover, the identity $\Im F_\psi (r)=0$ implies that 
 $|s(\lambda)|=1$.
 
 If the functions $v_{1}(x)$ and $v_{0}(x)$ decay  super-exponentially, then the remainder $o(1)$ in \e{eq:HL} can be replaced by a more precise term $- 2^{-1/2} b(\lambda)e^{- kx}+ O(e^{-\varkappa x})$ where $\varkappa$ is an arbitrary number. Relation \e{eq:uni} entails that
 \begin{equation}
 s(\lambda)\overline{b(\lambda)}=b(\lambda) 
\label{eq:HL1}\end{equation} 
which plays the role of identity \e{eq:4di7}.


\section{Perturbation by a  boundary condition}
 

Here we discuss the Hamiltonian  $H=D^4$ in the space $L_2({\Bbb R}_+)$ with some self-adjoint boundary conditions at the point $x=0$ and calculate   explicitly its resolvent.  In contrast to previous sections, we avoid here references to general results of scattering theory and   give direct proofs of all assertions.

 \medskip 

{\bf 5.1.}  
Self-adjoint extensions of a symmetric
 operator $D^4$  with domain $C_{0}^\infty ({\Bbb R}_{+})$ in the space $L_{2}({\Bbb R}_{+})$ are  defined 
by the formula $( Hu) (x)= u^{(4)} (x)$ on functions $u(x)$
from the Sobolev space
${\HH}^4({\Bbb R}_+)$ satisfying   some boundary   conditions at the point $x=0$. Let us describe all of them.   ``Generic" self-adjoint  boundary  conditions have the form
\begin{equation}\left\{\begin{array}{lcl}
u^{\prime\prime}(0)&=&\alpha u (0)+\alpha_1 u^{\prime }(0)
\\ 
u^{\prime\prime\prime}(0)&=&-\alpha_2 u (0) -\bar{\alpha} u^{\prime }(0),
\end{array}\right.
\label{eq:4ord}\end{equation}
where $\alpha_1$ and $\alpha_2$ are arbitrary real numbers and $\alpha$ is an arbitrary complex number. This family depends on four real constants. We introduce also ``exceptional"  
three-parameters 
  \begin{equation}
u'(0) =\alpha u(0), \q -u'''(0)+\bar{\alpha} u''(0) =\alpha_{2} u(0)  , 
\q \alpha \in{\Bbb C}, \q \alpha_{2}\in{\Bbb R}, 
 \label{eq:par2}\end{equation}
  and  a one-parameter 
 \begin{equation}
u(0) =0, \q u''(0)=\alpha_{1}u'(0) , \q\alpha_{1}\in{\Bbb R},
 \label{eq:par1}\end{equation}
 families of boundary conditions. To exhaust all self-adjoint extensions, we have to add      boundary condition \e{eq:4ordfr}. The operator $D^4$ with this boundary condition will be denoted $H_{0}$. We use also a special notation $H_{00}$ for the operator $D^4$ with
 boundary condition \e{eq:par1} where $\alpha_{1}=0$, that is 
 \begin{equation}
u(0) =0, \q u''(0)=0 .
 \label{eq:free}\end{equation}

  The quadratic form of the operator $H=H(\alpha,\alpha_1,\alpha_2)$ with boundary conditions \e{eq:4ord} is given by the expression
\begin{equation}
h[u,u] = \int_{0}^\infty |u'' (x)|^2 dx + \alpha_{2}|u(0) |^2+ 2\Re (\alpha  u(0)\bar{u}'(0)) 
+ \alpha_1 |u'(0) |^2. 
\label{eq:QFBC}\end{equation}
It is closed on the set ${\HH}^2 ({\Bbb R}_+)$. 
 In   case \e{eq:par2} the form is defined  by the expression
\begin{equation}
h[u,u] = \int_{0}^\infty |u'' (x)|^2 dx + \alpha_{2}|u(0) |^2  
\label{eq:QFBC1}\end{equation}
 on functions from $  {\HH}^2 ({\Bbb R}_+)$ satisfying the condition $u'(0)=\alpha u(0)$.  In   case \e{eq:par1}   expression \e{eq:QFBC} where $\alpha=\alpha_{2}=0$  remains true if the form is restricted on functions 
$u\in {\HH}^2 ({\Bbb R}_+)$ satisfying the condition $u(0)=0$.
Finally,  the quadratic  form of the operator $H_{0}$ is given by expression \e{eq:QFBC1} where $\alpha_{2}=0$ on functions $u\in {\HH}^2 ({\Bbb R}_+)$  such that $u(0)=u'(0)=0$.

Note that $ H \leq H_{0}$ for all boundary conditions because the quadratic form of the operator $H_{0}$  is defined on the smallest possible set. Similarly,   $H\leq H_{00}$ for boundary conditions \e{eq:4ord} and all $\alpha_{1}, \alpha_{2}, \alpha$ because the quadratic form of the operator $H_{00}$  is defined on a smaller (the boundary condition $u(0)=0$ is added) set than that of $H $.

 As a ``free" operator,  it is natural to
take the operator $H_0$ (the same operator as in  subs.~4.4). However technically it is more convenient to work with the  ``intermediary  free" operator $H_{00}$. 
 The reason for this choice of  the  free  operator is that $H_{00}$
 is the square of the operator $D^2$ with the boundary condition $u(0)=0$.

We   start with a construction of the
resolvent $R(z)$ of the operator $H$. As  a preliminary step, we shall find an explicit expression (cf. Lemma~\ref{2.1})
for the resolvent $R_{00}(z)$ of the operator $H_{00}$.

\begin{lemma}\label{4ord} 
Let $z\in{\Bbb C}\setminus [0,\infty)$ and $\zeta^4=z$, $\arg \zeta\in (0,\pi/2)$. Then
\begin{equation}
 (R_{00}(z)f)(x)=\frac{1}{4\zeta^3} \int_0^\infty(i e^{i\zeta |x-y|} -e^{-\zeta |x-y|}-i e^{i\zeta(x+y)}+  e^{-\zeta(x+y)})f(y)dy.
\label{eq:4ordr}\end{equation}
\end{lemma}

\begin{pf}
Using the Fourier sine transform, we see that
\[
 (R_{00}(z)f)(x)=2\pi^{-1}\int_0^\infty dk \sin(k x)(k^4-z)^{-1}
\int_0^\infty  dy \sin(k y) f(y) .
\]
Therefore the resolvent kernel equals
\[
  R_{00}(x,y;z) =-(2\pi)^{-1}\int_{-\infty}^\infty  
 e^{ik(x+y)}(k^4-z)^{-1}dk
+(2\pi)^{-1}\int_{-\infty}^\infty e^{ik | x-y|}(k^4-z)^{-1}dk.
\]
Both integrals can be complemented in the upper half-plane by a big half-circle and then calculated by residues at the points $\zeta$ and $i\zeta$.
\end{pf}

The solutions from $L_2({\Bbb R}_+)$ of the equation 
\[
u^{(4)}(x) =zu(x) + f(x)
\]
 satisfying different boundary conditions at the point $x=0$ differ by a linear combination of the solutions $e^{i\zeta x}$ and $e^{-\zeta x}$ of the homogeneous equation $u^{(4)}(x) =zu(x)$. Therefore we seek  the resolvent $R(z)$ of $H$ in the form (known as the Kre\u{\i}n formula)
\begin{equation} 
R(z)=R_{00}(z)+P(z),
\label{eq:4ord1}\end{equation}
where $P(z)$ is a two-dimensional operator defined by  equalities
\begin{eqnarray} 
(P(z)f)(x)=(p_{11}(\zeta) e^{i\zeta x}+p_{12}(\zeta) e^{-\zeta x} )Q_+(z)f
\nonumber\\
+ ( p_{21}(\zeta) e^{i\zeta x} +p_{22}(\zeta) e^{-\zeta x})Q_-(z)f
\label{eq:4ord2}\end{eqnarray}
and
\begin{equation}  
Q_+(z) f =\int_0^\infty e^{i\zeta y} f(y)dy, \quad
Q_-(z) f =\int_0^\infty e^{-\zeta y} f(y)dy.
\label{eq:4ord3}\end{equation}
Calculating   derivatives of expressions (\ref{eq:4ordr}) and (\ref{eq:4ord2}) at $x=0$, we find that
\begin{equation}   
\left\{\begin{array}{lcl} 
 (R(z)f)(0)&=& ( p_{11}(\zeta)  +p_{12}(\zeta))Q_+(z)f+( p_{21}(\zeta)  +p_{22}(\zeta))Q_-(z)f,
\\ 
(R(z)f)^\prime (0)&=& 2^{-1} \zeta^{-2}(Q_+(z)f-Q_-(z)f)
\\
&+&\zeta(( i p_{11}(\zeta)  -p_{12}(\zeta))
Q_+(z)f+(i p_{21}(\zeta) -p_{22}(z))Q_-(\zeta)f),
\\
 (R(z)f)^{\prime\prime}(0)&=& \zeta^2((- p_{11}(\zeta)  +p_{12}(\zeta))Q_+(z) f
 \\
 &+&(- p_{21}(\zeta) +p_{22}(\zeta))Q_-(z)f),
\\
(R(z)f)^{\prime\prime\prime} (0)&=& - 2^{-1}  (Q_+(z)f+Q_-(z)f)
\\
&-&\zeta^3(( i p_{11}(\zeta) 
+p_{12}(\zeta)) Q_+(z)f+(i p_{21}(\zeta) +p_{22}(\zeta))Q_-(z)f).
\end{array}\right.
\label{eq:RRx}\end{equation}

The coefficients $p_{jl}(\zeta)$ in \e{eq:4ord2} are determined by boundary conditions on the functions $u(x)= (R(z)f)(x)$ at the point $x=0$. 
We consider only  generic  boundary conditions (\ref{eq:4ord})    although formulas obtained above allow one to treat  also easily cases \e{eq:par2} and
  \e{eq:par1}.
Plugging   expressions \e{eq:RRx} into    (\ref{eq:4ord}) and equating coefficients at $Q_+ f$ and $Q_- f$, we obtain two systems of equations for $p_{11}$, $p_{12}$ and for
 $p_{21}$, $p_{22}$:
\begin{equation}\left\{\begin{array}{lcl}    
a_{11} (\zeta) p_{j1} (\zeta) + a_{12} (\zeta)p_{j2}(\zeta)=q_{j1}(\zeta),
\\
a_{21}(\zeta)p_{j1} (\zeta) +a_{22} (\zeta)p_{j2} (\zeta)=q_{j2} (\zeta),
\end{array}\right.
\label{eq:4ord5}\end{equation}
where $j=1,2$. The coefficients of these systems are given by explicit formulas 
\begin{equation}\left\{\begin{array}{lcl}    
a_{11} (\zeta)=\alpha+i\alpha_1\zeta +\zeta^2,
\quad
a_{12} (\zeta)=\alpha -\alpha_1\zeta -\zeta^2,
\\ 
a_{21} (\zeta)=-\alpha_2 -i\bar{\alpha}\zeta +i\zeta^3,
\quad  
a_{22} (\zeta)= -\alpha_2 +\bar{\alpha}\zeta + \zeta^3,
\end{array}\right.
\label{eq:4ord6}\end{equation}
and
\[
\left\{\begin{array}{lcl}    
 q_{11}(\zeta)= -(2\zeta^2)^{-1}\alpha_1,\quad
 q_{12}(\zeta)=  (2\zeta^2)^{-1}\bar{\alpha}-2^{-1},
\\
q_{21}(\zeta)= (2\zeta^2)^{-1}\alpha_1,\quad
 q_{22}(\zeta)= - (2\zeta^2)^{-1}\bar{\alpha}-2^{-1}.
\end{array}\right.
\]
Of course,   system (\ref{eq:4ord5}) can   easily be solved. Let us set
$\alpha_0= \alpha_1 \alpha_2-|\alpha|^2$ and
\begin{eqnarray}
\Omega(\zeta)&=& -2^{-1/2}   e^{-\pi i/4}  \zeta^{-1}(a_{11}(\zeta)a_{22}(\zeta)- a_{12}(\zeta) a_{21}(\zeta)) 
\nonumber\\
&=&  \alpha_0  + (1-i) \alpha_2\zeta  + 2i\Re \alpha\,
 \zeta^2 -(1+i) \alpha_1\zeta^3      -  \zeta^4.
\label{eq:4ord8}\end{eqnarray}
Then
\begin{equation}
p_{jl}(\zeta)=  2^{-3/2} e^{-\pi i/4} \zeta^{-3} \Omega (\zeta)^{-1}  {\bf p}_{jl}(\zeta),\q j,l=1,2,
\label{eq:P}\end{equation}
where 
\begin{equation}   
\left\{\begin{array}{lcl}    
  {\bf p}_{11}(\zeta)=     
-  \alpha_0 -2\Re \alpha\,  \zeta^2+2\alpha_1 \zeta^3+\zeta^4,
 \\
  {\bf p}_{12}(\zeta)=    \alpha_0 + 2i\Im \alpha 
\,\zeta^2 +\zeta^4,
 \\
  {\bf p}_{21}(\zeta)=  \alpha_0 - 2i\Im \alpha\, 
\zeta^2 +\zeta^4,
 \\
  {\bf p}_{22}(\zeta)=   -  \alpha_0 + 2\Re \alpha \,
\zeta^2 + 2 i \alpha_1 \zeta^3+\zeta^4.
\end{array}\right.
\label{eq:4ord9}\end{equation}
Note that
\[
\overline{\Omega(\zeta )}=\Omega(i\bar{\zeta }), \q  \zeta \in   {\Bbb C },
\]
so that zeros of the function $\Omega(\zeta) $  are symmetric with respect to the line $\Re \zeta=\Im \zeta$.

Let us summarize the results obtained.

\begin{theorem}\label{4ord1}
  Let $z\in{\Bbb C}\setminus [0,\infty)$ and $\zeta^4=z$, $\arg
\zeta\in (0,\pi/2)$. Then the resolvent $R(z)$ of the 
operator $H=d^4/dx^4$ with  boundary conditions $(\ref{eq:4ord})$, \e{eq:par2}
or \e{eq:par1} is given by formula
$(\ref{eq:4ord1})$ where the resolvent $R_{00}(z)$     of the 
operator $H_{00}$ is determined by formula
$(\ref{eq:4ordr})$ and the operator $P(z)$  is determined by formulas
$(\ref{eq:4ord2})$ and $(\ref{eq:4ord3})$. In   case $(\ref{eq:4ord})$ the coefficients $p_{jl}(\zeta)$ are defined by formulas
    $(\ref{eq:4ord8})$, $(\ref{eq:P})$ and $(\ref{eq:4ord9})$.
\end{theorem}

Note a particular case of this result. If $\alpha=\alpha_{1}=\alpha_{2}=0$, then $\Omega(\zeta)=-\zeta^4$ and 
 \begin{equation}
P(x,y;z)= -2^{-3/2} e^{-\pi i /4} \zeta ^{-3}  (e^{i\zeta x}+ e^{-\zeta x})(e^{i\zeta y}+e^{-\zeta y}).
 \label{eq:000}\end{equation}

Of course explicit formulas for the coefficients $p_{jl}(\zeta)$ in \e{eq:4ord2} can easily be written down also for other boundary conditions. For example, let us consider the operators  $H_0$ and $H_{1}$ determined by boundary conditions   \e{eq:4ordfr} and
  \begin{equation}
  u'(0)=u'''(0)=0,
 \label{eq:ex2}\end{equation}
 respectively. Using formulas 
 \e{eq:RRx}, we see that their resolvents are given  by the
 equality 
 \begin{equation}
R_{j}(z)= R_{00}(z)+P_{j}(z), \q j=0,1,
 \label{eq:exP1}\end{equation}
   where $P_j(z)$ are integral operators  with kernels 
 \begin{equation}
P_0(x,y;z)= (4\zeta^3)^{-1}(i-1) (e^{i\zeta x}-e^{-\zeta x})(e^{i\zeta y}-e^{-\zeta y})
 \label{eq:exP1a}\end{equation}
 and 
  \begin{equation}
P_1(x,y;z)= (2\zeta^3)^{-1}  ( i e^{i\zeta (x+y)}-e^{-\zeta (x+y)}).
 \label{eq:exP2}\end{equation}

 \medskip

{\bf 5.2.}
Let us discuss eigenvalues of the operator $H$ corresponding to  boundary conditions  (\ref{eq:4ord}). Recall that the function $\Omega$ was defined by formula \e{eq:4ord8}.

 \begin{proposition}\label{negX}
 Zeros of the function $\Omega(\zeta)$ from the sector $\arg \zeta\in (0,\pi/2)$ lie on the ray  $\arg \zeta =\pi/4$. A point $\lambda=-k^4$, $k>0$,   is an eigenvalue of the operator $H$ if and only if    $\Omega( e^{i\pi/4} k)=0$, that is
 \[
k^4+\sqrt{2}\alpha_{1}k^3 -2\Re \alpha \, k^2 + \sqrt{2} \alpha_{2} \, k + \alpha_{0}=0.
\] 
    Multiplicities of an eigenvalue $\lambda$ of the operator $H$  and of the   zero $e^{i\pi/4} k$ of the function $\Omega $ coincide.  A point
 $\lambda_{0}=-k^4_{0}$ is a degenerate $($of multiplicity $2)$ eigenvalue of the operator $H$ if and only if  
   \begin{equation}
 \alpha_{1}= -\sqrt{2}k_{0}, \q \alpha=-k_{0}^2, \q  \alpha_{2}= - \sqrt{2}k_{0}^3.
\label{eq:mult1}\end{equation}
 In this case the eigenfunctions are defined by equations
$$
\psi_{1}(x)=\exp ((-1+i) k_{0} x/ \sqrt{2}), \q \psi_{2}(x)=\exp ((-1- i) k_{0} x/ \sqrt{2}) 
$$
and necessarily
   \begin{equation}
\Omega(\zeta)= -(\zeta -e^{i\pi/4} k_{0})^2  (\zeta^2+ i k_{0}^2).
\label{eq:mult2}\end{equation}
 \end{proposition} 

 \begin{pf}
   Clearly,  a solution
 \begin{equation}
u(x)= p_{1} e^{i\zeta x} +  p_{2} e^{- \zeta x}
 \label{eq:4or5}\end{equation}
of   the equation
 \begin{equation}
  u^{(4)}(x)= \zeta^4 u (x)
   \label{eq:4or5h}\end{equation}
   satisfies boundary conditions (\ref{eq:4ord}) if and only if (cf. \e{eq:4ord5})
 \begin{equation}\left\{\begin{array}{lcl}    
a_{11} (\zeta) p_{1}   + a_{12} (\zeta)p_{2} = 0,
\\
a_{21}(\zeta)p_{1}   + a_{22} (\zeta)p_{2}  =0,
\end{array}\right.
\label{eq:4ord5h}\end{equation}
where the coefficients $a_{jl}(\zeta)$ are defined by formulas \e{eq:4ord6}.
Since both functions $e^{i\zeta x} $ and $e^{- \zeta x}$ belong to $L_{2}({\Bbb R}_{+})$,
the point  $\lambda =\zeta^4$ is an eigenvalue of the operator $H$ if and only if system
\e{eq:4ord5h} has a non-trivial solution, that is  $\Omega(\zeta)=0$. It follows that   $\arg \zeta = \pi/ 4$ because necessarily  $ \lambda< 0$.
 
 A point $\lambda_{0}=-k^4_{0}$ is an    eigenvalue  of multiplicity 2 if and only if system \e{eq:4ord5h}  where $\zeta_{0}=e^{\pi i/4} k_{0}$ is satisfied for all numbers $p_{1} $ and $p_{2}$, that is
  \begin{equation}
  a_{11} (\zeta_{0})=a_{12} (\zeta_{0})=0, \q a_{21} (\zeta_{0})=a_{22} (\zeta_{0})=0. 
  \label{eq:mult5}\end{equation}
 Using \e{eq:4ord6} and solving the first two    equations \e{eq:mult5}, we obtain expressions
 \e{eq:mult1} for $ \alpha_{1}$ and $ \alpha$. Similarly, considering the last two  equations \e{eq:mult5}, we obtain expression
 \e{eq:mult1} for $ \alpha_{2}$. Plugging  expressions  \e{eq:mult1} into \e{eq:4ord8}, we get representation \e{eq:mult2} so that the point $\zeta_{0}$ is a double zero of the function $\Omega(\zeta)$. 
 
 Conversely, if $\Omega(\zeta_{0})=\Omega'(\zeta_{0})=0$, then ${\bf p}_{jl}(\zeta_{0})=0$ for all $j,l=1,2$ because the functions $p_{jl}(\zeta )$ can have only simple poles. This is a consequence of representations \e{eq:4ord1}, \e{eq:4ord2} for the resolvent $R(z)$ whose  poles are simple. By virtue of \e{eq:4ord9} the equations ${\bf p}_{12}(\zeta_{0})={\bf p}_{21}(\zeta_{0})= 0$ yield $\Im \alpha=0$, $\alpha_0    +\zeta^4_{0}=0$.
Therefore  the equations ${\bf p}_{11}(\zeta_{0})={\bf p}_{22}(\zeta_{0})= 0$ can be written as
\[  
\left\{\begin{array}{lcl}    
-    \alpha  +\alpha_1 \zeta_{0}+\zeta^2_{0} = 0,
 \\
    \alpha  + i \alpha_1 \zeta_{0}+ \zeta^2_{0} =  0,
\end{array}\right.
\]
whence  $(1+i) \alpha_{1}+2\zeta_{0}=0$ and  $(1- i) \alpha_{1}\zeta_{0}=2 \alpha$.
 This implies equations \e{eq:mult1}   so that equations \e{eq:mult5} for numbers \e{eq:4ord6} are also satisfied. Thus,  function \e{eq:4or5} where $\zeta=\zeta _{0}$ satisfies boundary conditions \e{eq:4ord} for all $p_{1}$, $p_{2}$ and hence $-k_{0}^4$ is an eigenvalue of multiplicity $2$.
  \end{pf}

 Since the rank of the operator $R(z)-R_{00}(z)$ equals $2$ and $H_{00}\geq 0$, the operator $H$ might have at most $2$ negative eigenvalues with multiplicity taken into account. The following result (cf. \cite{Y1}) makes this assertion more precise.
 
 \begin{proposition}\label{neg}
The total numbers of negative eigenvalues $($counted with their  multiplicity$)$ of the operators $H$ and
\[
 A= \left(\begin{array}{cc}
 \alpha_{2}&\bar{\alpha}
\\ \alpha &\alpha_{1}
\end{array} \right): {\Bbb C}^2\to  {\Bbb C}^2
\]
   coincide.
\end{proposition} 
 
  \begin{pf}
  Let $E({\Bbb R}_{-})$ and ${\cal E}({\Bbb R}_{-})$ be spectral projectors of the operators $H$ and $A$, respectively, corresponding to the set ${\Bbb R}_{-}=(-\infty, 0)$.
  We have to check that
  \[
  N:= \dim E({\Bbb R}_{-}){\cal H}= \dim {\cal E}({\Bbb R}_{-}) {\Bbb C}^2=:n.
  \]
  Let us define the mapping $J:{\HH}^2 ({\Bbb R}_{+})\to {\Bbb C}^2$ by the relation
  $ J u=(u(0), u'(0))$. Then formula \e{eq:QFBC} can be written as
  \begin{equation} 
h[u,u] = \int_{0}^\infty |u'' (x)|^2 dx + \langle AJu,Ju\rangle.
\label{eq:neg1}\end{equation}
  It follows that if $u\in E({\Bbb R}_{-}){\cal H}$, $u\neq 0$, then $\langle AJu,Ju\rangle < 0$. Thus, the quadratic form of the operator $A$ is negative on the subspace $JE({\Bbb R}_{-}){\cal H}$.  This subspace has dimension $N$ because $Ju\neq 0$ for $u\neq 0$ and hence $n\geq N$.

Conversely, pick a function $\varphi\in C^\infty ({\Bbb R}_{+}) $ such that $\varphi (x)=1$ in a neighborhood of the point $x=0$ and $\varphi (x)=0$ for sufficiently large $x$ and set
\[
(\tilde{J}_{\varepsilon}a) (x)= (a_{1}+ a_{2}x) \varphi(\varepsilon x), \q a=(a_{1}, a_{2}), \q \tilde{J}_{\varepsilon}:   {\Bbb C}^2 \to {\HH}^2 ({\Bbb R}_{+}).
\]
Remark that there exists a	constant $\gamma>0$ such that 
$\langle A a, a\rangle \leq -\gamma {\pmb |}  a {\pmb |} ^2 $ for  all $a \in {\cal E}({\Bbb R}_{-}) {\Bbb C}^2$. According to \e{eq:neg1} we have that
  \begin{equation} 
h[\tilde{J}_{\varepsilon}a ,\tilde{J}_{\varepsilon}a] = \int_{0}^\infty |(\tilde{J}_{\varepsilon}a)'' (x)|^2 dx + \langle A a, a\rangle,
\label{eq:neg2}\end{equation}
where 
\[
\int_{0}^\infty |(\tilde{J}_{\varepsilon}a)'' (x)|^2 dx
 \leq \varepsilon {\pmb |}  a {\pmb |} ^2.
\]
Thus,   expression \e{eq:neg2} is negative if   $\varepsilon<\gamma$ so that 
$\tilde{J}_{\varepsilon} {\cal E}({\Bbb R}_{-}) {\Bbb C}^2$ is a
subspace of dimension $n$ on which the quadratic form $h$ is negative. 
This implies that $n\leq N$. 
\end{pf}

\begin{corollary}\label{neg1}
If $\alpha_{0}=\alpha_{1}\alpha_{2}-|\alpha|^2<0$, then the operator  $H$ has precisely one negative eigenvalue.
If $\alpha_{0}>0$, then the operator  $H$ has two negative eigenvalues for  $\alpha_{1}<0$ $($or equivalently $\alpha_{2}<0)$ and   $H\geq 0$  for $\alpha_{1}>0$ $($or equivalently $\alpha_{2}>0)$.
\end{corollary} 

\begin{pf}
If $\alpha_{0}=\det A < 0$, then eigenvalues $\mu_{1}$ and $\mu_{2}$ of the operator $A$ have different signs. If $\alpha_{0}=\det A > 0$, then   $\mu_{1} \mu_{2}>0 $, $\alpha_{1} \alpha_{2}>0 $ and $\mu_{1} + \mu_{2}=\alpha_{1}+ \alpha_{2}  $. Thus, all  four numbers  
$\mu_{1}, \mu_{2}   $ and $\alpha_{1}, \alpha_{2}  $ have the same sign.
\end{pf}
 
 Let us further  consider  zeros of the function $\Omega(\zeta)$ on the half-axis $\zeta =k$ (or equivalently  $\zeta =ik$) where $k>0$.  
   
\begin{proposition}\label{4ord1q}
The function $\Omega(k)=0$ for $k>0$ if and only if  $ \lambda=k^4 $
   is an eigenvalue of the operator $H$.
\end{proposition}  

\begin{pf}
Put $\zeta=k$.
Recall that   solution \e{eq:4or5} of equation \e{eq:4or5h}
satisfies boundary conditions
 (\ref{eq:4ord}) if and only if the coefficients $p_{1}$ and $p_{2}$ satisfy system
\e{eq:4ord5h}.
If $\lambda $ is an eigenvalue of  $H$, then   $\psi (x)=\exp (- k x)$ is the eigenfunction
of $H$ and hence system \e{eq:4ord5h} is satisfied with $p_{1}=0$ and $p_{2}=1$. It follows that $a_{12}(k)= a_{22}(k)=0$ and hence $\Omega(k)=0$.

 Conversely,     if $\Omega(k)=0$, then   system \e{eq:4ord5h}
has a nontrivial solution $ p_{1} ,  p_{2} $ so that solution \e{eq:4or5} of equation \e{eq:4or5h} satisfies boundary conditions
 (\ref{eq:4ord}). It remains only to show that $ p_{1} =0 $.
 Let us  multiply (cf. the proof of Lemma~\ref{4di})  equation \e{eq:4or5h} by $ \overline{u(x)}$ and integrate it over an interval $(0, r)$. Then we integrate by parts and take the imaginary part. Since the non-integral terms    disappear at $x=0$, we obtain the identity
\[
\Im (u'''(r)\overline{u(r)} - u''(r)\overline{u'(r)})=0.
\]
Applying it to function \e{eq:4or5} and neglecting terms exponentially decaying as $r\to\infty$, we see that $p_{1} =0$ and hence  $u(x)=e^{-kx}$ is an eigenfunction  of the operator $H$.
\end{pf}
  
  Positive eigenvalues are not exceptional for the operator $H$.  It is easy to give simple necessary and sufficient conditions for their existence.   
  
  \begin{proposition}\label{4ord1qq}
A point  $ \lambda =k^4  $, $ k>0$,    is an eigenvalue of the operator $H$ if and only if 
\begin{equation}
\alpha=\overline{\alpha}\q \mathrm{and}\q
 \alpha_{1}= (\alpha-k^2) k^{-1} , \q \alpha_{2}=  (\alpha+ k^2 ) k .
  \label{eq:4or5a}\end{equation}  
  In this case
   \begin{equation}
  \Omega'(k ) =-2 k  (\alpha-ik ^2)\neq 0 
 \label{eq:4tr16}\end{equation}
 so that $k $ is  a simple zero of the function $\Omega(\zeta)$.
\end{proposition}

 \begin{pf}
 The eigenfunction of $H$  is necessarily $e^{-k x}$, and boundary conditions  
 \e{eq:4ord}  for this function are equivalent to   equations    \e{eq:4or5a}. Differentiating \e{eq:4ord8} and using expressions \e{eq:4or5a} for $\alpha_{1}$ and $\alpha_{2}$, we arrive at \e{eq:4tr16}.
  \end{pf}
  
    \begin{corollary}\label{4ord1qqw}
The operator $H$ cannot have more than one positive eigenvalue. 
\end{corollary}

 \begin{pf}
 If $\lambda=k^4$, $k>0$, is an eigenvalue of $H$, then it follows from   equations \e{eq:4or5a} that $k^4+\alpha_{0}=0$ which determines $k>0$ uniquely.
  \end{pf}

  Thus, for each given $\lambda >0$, there is a one-dimensional manifold in the 
  four-dimensional space of parameters $\alpha  , \alpha_{1}, \alpha_{2}$ parametrized by $\alpha=\overline{\alpha}$ such that the corresponding operators $H=H(\alpha,\alpha_1,\alpha_2)$ have  an eigenvalue at the point $\lambda $. If $k$ varies over 
  ${\Bbb R}_{+}$, then   equations \e{eq:4or5a}  determine (parametrically) a surface in the space of parameters  such that the operators $H$ have a positive eigenvalue. Note that the condition $\alpha_{0}<0$ is necessary for the existence of a positive eigenvalue.

 \medskip

{\bf 5.3.}
Next  we   calculate the spectral measure $E(\lambda)$ of the operator $H$ corresponding to boundary conditions \e{eq:4ord}. 

\begin{proposition}\label{4ord3}
Suppose that  $\lambda=k^4>0$ is not an eigenvalue of the operator $H$. Set
\begin{equation}
  s(\lambda)= \overline{\Omega(k)} /  \Omega(k) 
\label{eq:4ord17}\end{equation}
and
\begin{equation}
b(\lambda)=(\alpha_0+2i\Im\alpha \, k^2+k^4)/ \Omega(k) 
\label{eq:4ord17za}\end{equation}
where the function $\Omega(k)$ is defined by formula $(\ref{eq:4ord8})$. Then
$ dE(\lambda)/d\lambda$ is the integral operator with kernel
\begin{equation}
  dE(x,y;\lambda)/d\lambda=(2\pi )^{-1} k^{-3} \psi(x,\lambda)\overline{\psi(y,\lambda)},
\label{eq:4ord14}\end{equation}
where 
\begin{equation}
  \psi(x,\lambda)=2^{-1}(s(\lambda)e^{ik x+\pi i/4}+ e^{-ik x-\pi i/4})
  - 2^{-1/2} b(\lambda) e^{-k x}.
\label{eq:4ord15}\end{equation}
\end{proposition}

\begin{pf}
Passing to the limit $z\rightarrow\lambda\pm i0$ in equality (\ref{eq:4ord1})  and using the relation between boundary values of the resolvent and the spectral measure (see formula \e{eq:Cauchy}), we  find  that
\begin{equation}
dE(x,y; \lambda)/d\lambda= dE_{00}(x,y; \lambda)/d\lambda
+  (2 \pi  i)^{-1} ( P (x,y; \lambda + i0)-\overline{P (y,x; \lambda + i0)}).
\label{eq:4ord11}\end{equation}
The kernel $P (x,y; z)$ of the  operator $P (z)$ is defined by formulas (\ref{eq:4ord2}) and (\ref{eq:4ord3}).
Since $E_{00} (\lambda)=E_{\sqrt{H_{00}} }(\sqrt{\lambda} )$ where $\sqrt{H_{00}} =
D^2$ with the boundary
condition $u(0)=0$, we have that  
\[ 
dE_{00}(x,y; \lambda)/d\lambda=(2\pi )^{-1} k^{-3} \sin k x \sin k y.
\]
 Therefore it follows from relation \e{eq:4ord11} that 
   \begin{eqnarray}
 dE(x,y;\lambda)/d\lambda= -(8 \pi )^{-1} k^{-3} (e^{ik x}-e^{-ik x})(e^{ik y}-e^{-ik y})
\nonumber\\ 
+ (2\pi i)^{-1}\Big( p_{11}(k) e^{ik (x + y)}  - \overline{p_{11}(k)}  e^{-ik (x + y)}
\nonumber\\ 
 + p_{12}(k) e^{- k (x - iy)}
- \overline{p_{21}(k)}     e^{ -k (x+i y )}             
    + p_{21}(k) e^{- k (-ix +y)} - \overline{p_{12}(k) }  e^{-k (i x+ y )} 
  \nonumber\\
 + ( p_{22}(k) -\overline{p_{22}  (k)}  ) e^{-k (x + y)}\Big),
  \label{eq:4ord13}\end{eqnarray}
where   the coefficients $p_{jl}$ are defined by formulas
    (\ref{eq:4ord8}) -- (\ref{eq:4ord9}).
    
    Let us plug  \e{eq:4ord15} into the right-hand side of  \e{eq:4ord14}. Taking into account definitions  \e{eq:4ord17} and  \e{eq:4ord17za}, we see that the coefficients at all terms $e^{ik (x + y)}$, $e^{-ik (x + y)}$, $e^{ik (-x + y)}$, $e^{ik (x - y)}$,  $e^{- k (x - iy)}$, $ e^{ -k (x+i y )}$, $e^{- k (-ix +y)}$,  $e^{- k (ix +y)}$
    and $e^{- k (x + y)}$ are the same as in \e{eq:4ord13}. This proves relation \e{eq:4ord14}.
    \end{pf}
    
    The functions $\psi  (x,\lambda) $ defined by formula \e{eq:4ord15} satisfy  of course  boundary conditions \e{eq:4ord}. They are known as eigenfunctions of the continuous spectrum of the operator $H$.     These functions describe an interaction of a plane wave $e^{-ikx-\pi i/4}$ coming from $+\infty$ with a ``point" potential at $x=0$. The coefficient $s(\lambda)$ at the reflected wave $e^{ikx+\pi i/4}$ is known as the scattering matrix.
    Note that if $\Im\alpha=0$ (in this case the problem is not only self-adjoint but is also real), then    according to equalities (\ref{eq:4ord17}) and (\ref{eq:4ord17za}) the functions $s(\lambda) $ and $b(\lambda)$ are related by formula \e{eq:HL1}. Thus,  the scattering matrix $s(\lambda) $  can be recovered from the coefficient $b(\lambda)$ at the exponentially decaying mode of eigenfunction (\ref{eq:4ord15}). 
    
    Proposition~\ref{4ord3} implies that the positive spectrum of the operator $H$ is absolutely continuous except, possibly, a single eigenvalue  $ \lambda_{0}$. Moreover, we have the following result.
    
    \begin{proposition}\label{4ord5sp}
     The functions  $s(\lambda) $ and  $b(\lambda) $  are infinitely differentiable 
 for all $\lambda>0$. If  $\lambda_{0} $ is a positive eigenvalue of the operator $H$,
then
\begin{equation}
s(\lambda_{0}) =  \frac{ \alpha+ i\sqrt{\lambda_{0}} } { \alpha-i\sqrt{\lambda_{0}} }, \q
b(\lambda_{0}) = -2 \frac{ \sqrt{\lambda_{0}} } { \alpha-i\sqrt{\lambda_{0}} }.
 \label{eq:4trsm}\end{equation}
  \end{proposition}
  
  \begin{pf}
  Let us proceed 
 from formulas (\ref{eq:4ord17}) and \e{eq:4ord17za}. Clearly,  the functions  $s(\lambda) $ and  $b(\lambda) $  are infinitely differentiable 
 away from the point $ \lambda_{0}=k_{0 }^4$ (the eigenvalue of $H$) where $\Omega(k_{0})=0$. As far as their behavior at $\lambda_{0}$ is concerned, the result about $s(\lambda)$ follows from formulas \e{eq:4tr16} and
  \[
  s(\lambda_{0}) =\overline{ \Omega'(k_{0}) }/\Omega'(k_{0}).
  \]
  To consider $b(\lambda)$, we take additionally into account  that according to
 Proposition~\ref{4ord1qq} a point  $ \lambda_{0}=k^4_{0} $   is a positive eigenvalue of the operator $H$ if and only if   conditions \e{eq:4or5a} are satisfied for $k=k_{0}$.  Therefore the numerator in \e{eq:4ord17za}  equals zero at $k=k_{0}$ and its derivative equals  $4 k_{0}^3$. This yields relation \e{eq:4trsm} for $b(\lambda_{0})$.
   \end{pf}
   
   We emphasize that both functions $i s(\lambda_{0}) e^{i k_{0}x} +  e^{- i k_{0}x}$ and 
   $e^{- k_{0}x}$ satisfy boundary condition \e{eq:4ord} if $\lambda_{0}$ is an eigenvalue.

    \begin{remark}\label{RemB}
     The function    $b(\lambda) $  might have zero only at one point $\lambda_{0}=-\alpha_{0}$.  According to Proposition~\ref{4ord1qq} this condition is satisfied if $\lambda_{0}$ is an eigenvalue of $H$. However according to   Proposition~\ref{4ord5sp} in this case $b(\lambda_{0})\neq 0 $. On the contrary, $b(\lambda_{0}) = 0 $ if $\lambda_{0}=-\alpha_{0}$, $\alpha=\bar{\alpha}$, but the last two conditions \e{eq:4or5a}  are violated at $k=k_{0}$. For example, if $\alpha=k_{0}^2$, $\alpha_{1}=0$ but $\alpha_{2} \neq 2 k_{0}$, then $b(\lambda_{0}) =0$  for $\lambda_{0}= k_{0}^4$.
       \end{remark}     
    
 \medskip

{\bf 5.4.}
Now   we are in a position to establish
 an expansion in eigenfunctions  of the operator $H$.  

\begin{theorem}\label{4ord4} 
 Let $\psi_-(x,\lambda ) = \psi(x, \lambda)$ be defined by formulas \e{eq:4ord17},  \e{eq:4ord17za},  \e{eq:4ord15}, and let 
 \begin{equation}
\psi_+ (x,\lambda) = \overline{s(\lambda)}  \psi_{-}(x,\lambda) . 
\label{eq:psi}\end{equation}
 Define on the set
$  L_1({\Bbb R}_+)\cap L_2 ({\Bbb R}_+)$ the mappings ${\cal F}_\pm$   by  the formula  
 \begin{equation}
({\cal F }_{\pm} f)(\lambda)= 2^{-1}(2\pi)^{-1/2} \lambda^{-3/8} \int_0 ^\infty \overline{\psi_{  \pm}(x,\lambda)}f(x) dx . 
\label{eq:speHL}\end{equation}
 These mappings extend to bounded operators on the space
 $ L_2({\Bbb R}_+)$ and satisfy relations \e{eq:F1}   and \e{eq:F2}.
\end{theorem}

\begin{pf}
Intertwining property \e{eq:F2} holds because the functions $\psi_\pm (x,\lambda) $ satisfy the equation $\psi_\pm^{(4)} =\lambda  \psi_\pm$  and boundary conditions \e{eq:4ord}.
The second relation  \e{eq:F1} is obtained by integration of  representation \e{eq:4ord14} over $\lambda\in{\Bbb R}_{+}$.  Now for the proof of   the first relation  \e{eq:F1}, we have to check that      the kernel of the operator ${\cal F}^*_{\pm}$ is trivial. Supposing that ${\cal F}^*_{\pm}g=0$ and hence ${\cal F}^*_{\pm}E_{A}(\alpha,\beta)g=0$ for all 
$ \alpha, \beta\in {\Bbb R}_{+}$, we have
\[
\int_\alpha^\beta \psi (x,\lambda) g(\lambda) d \lambda=0, \q \forall x\geq 0. 
\] 
 Let us differentiate twice this relation, use \e{eq:4ord15} and set $x=0$. Since $\alpha$ and $\beta$ are arbitrary, we obtain that 
 \begin{eqnarray*}
(s(\lambda)-i -2^{1/2} e^{-\pi i/4} b(\lambda) ) g(\lambda)&=& 0,\quad (s(\lambda)+i -2^{1/2} e^{\pi i/4} b(\lambda) ) g(\lambda)=0,
\\
 (s(\lambda)-i + 2^{1/2} e^{-\pi i/4} b(\lambda) )  g(\lambda)&=&0
\end{eqnarray*}
for a.e.  $\lambdaÊ>0$. If $ g(\lambda)\neq 0$ for some $\lambda$, then the first and third equations imply that 
$s(\lambda)=i$ and $b(\lambda)=0$  which contradicts the second equation.
Thus,   $ g(\lambda) =0$ for a.e.  $\lambda>0$.  
\end{pf}

 \begin{remark}\label{RemBb}
If $\alpha=\bar{\alpha}$, then $\psi_+ (x,\lambda) = \overline{\psi_{-}(x,\lambda)}  $.
       \end{remark}

Consider now the operator $H_{0}=D^4$ corresponding to boundary conditions
 \e{eq:4ordfr}. Its eigenfunctions are defined by   formula \e{eq:5ordfr}.
 It follows from formulas \e{eq:exP1} and \e{eq:exP1a}  that the spectral measure 
 $E_{0}(\lambda)$ of the operator $H_{0}$ satisfies relation \e{eq:4ord14} where the role of $\psi$ is played by $\psi_{0}$. Theorem~\ref{4ord4} applies of course to the operator ${\cal F}_{0}$ defined by the formula
  \begin{equation}
({\cal F }_{0} f)(\lambda)= 2^{-1}(2\pi)^{-1/2} \lambda^{-3/8} \int_0 ^\infty  \psi_{ 0}(x,\lambda) f(x) dx . 
\label{eq:speHLfr}\end{equation}
Moreover, the operator ${\cal F}_{0}$ is unitary because the operator $H_{0}$ does not have eigenvalues.

Theorem~\ref{4ord4}    allows us to
 construct directly the time-dependent scattering theory in the same way as for the second order differential operators. We need also the following auxiliary assertion.
 
 \begin{lemma}\label{invprA}
If $ u \in C_{0}^\infty ({\Bbb R}_{+})$,
then
\begin{equation}
 \lim_{t \rightarrow\pm\infty}\int_0^\infty d x
\,\Bigl|
\int_0^\infty\exp(\mp i k x-i k^4 t)u (k)dk\Bigr|^2=0
\label{eq:invpr}\end{equation}
and
\begin{equation}
 \lim_{| t |\rightarrow\infty}\int_0^\infty d x
\,\Bigl|
\int_0^\infty\exp(- k x-i k^4 t)u (k)dk\Bigr|^2=0.
\label{eq:invprA}\end{equation}
\end{lemma}

\begin{pf}
Both relations \e{eq:invpr} and \e{eq:invprA}
are obtained by  a direct integration by parts which shows that
 the  integral over $k$ is bounded by $C(x+|t|)^{-1}$. 
 \end{pf}

\begin{theorem}\label{4ord5} 
  The  wave operators  $W_\pm=W_\pm (H,H_{0})$ exist, are complete and satisfy the equality $W_\pm={\cal F}_{\pm}^* {\cal F}_{0}$.
The scattering operator ${\cal S}$ for the pair $H_{0}$, $H$ acts in the space $L_{2}({\Bbb R}_{+})$ as multiplication by the function $s(\lambda)$   defined by
formulas $(\ref{eq:4ord8})$, $(\ref{eq:4ord17})$.
\end{theorem}

\begin{pf}
 According to Theorem~\ref{4ord4} all results about the wave operators $W_\pm  $ follow from   the   relation 
 \begin{equation}
 \lim_{t \rightarrow\pm\infty}\| ({\cal F}_{\pm}^* -{\cal F}_{0}^*) \exp(-iAt) g\| =0
\label{eq:2.2.29}\end{equation}
where  $g (k) $ is an arbitrary function from $L_2({\Bbb R}_+)$. 
It suffices to check (\ref{eq:2.2.29}) on the  set 
$C_0^\infty({\Bbb R}_+ )$.
Using (\ref{eq:speHL}), (\ref{eq:speHLfr}), we obtain that
\begin{eqnarray}
 (({\cal F}_\pm^\ast - {\cal F}_{0}^\ast)\exp (-iAt)g) (x)
\nonumber\\  = 2^{-1} (2\pi)^{-1/2}\int_0^\infty  (\psi_\pm (x,\lambda)-\psi_0 (x,\lambda))\exp(-i \lambda  t)g (\lambda) \lambda^{-3/8} d\lambda.
\label{eq:2.2.30xx}\end{eqnarray}
 Let us check (\ref{eq:2.2.29}), for example, for the sign $``-"$. It follows from \e{eq:5ordfr} and (\ref{eq:4ord15})   that
\[
    \psi_-(x, \lambda)-\psi_0(x, \lambda)  = 2^{-1 }  (s(\lambda)-1)e^{ik x+\pi i/4}
     - 2^{-1/2}  ( b(\lambda)- 1) \exp(-k x).
\]
The contributions  to (\ref{eq:2.2.30xx}) of the first and
second  terms in the right-hand side   tend in $L_2({\Bbb R}_+)$  to zero as $t\ri-\infty$   according to relations \e{eq:invpr} and \e{eq:invprA}, respectively. 

Equality $W_\pm={\cal F}_{\pm}^* {\cal F}_{0}$ implies that the scattering operator ${\cal S}$ satisfies relation \e{eq:S}. It remains to remark that 
   $ ({\cal F}_+ f)(\lambda) =s(\lambda) ({\cal F}_- f)(\lambda)$ according to \e{eq:psi}.
\end{pf}

 \begin{remark}\label{invprAB}
 Since the operator ${\cal F}_{0}$ as well as  the Fourier transform are bounded operators, it follows from formula \e{eq:5ordfr} that the integral operator $T$ with kernel $\exp(- k x)$ is bounded in $L_2({\Bbb R}_+)$. This result is a by-product of our considerations, but it is not of course new. Indeed, we have that 
 \[
\| T u \|^2= ({\cal C} u , u)
\]
where ${\cal C}$ is the integral operator   with kernel $ ( k +k')^{-1}$ known as Carleman's operator. A proof of its boundedness can be found, e.g., in \cite{Pow}. 
We note also that boundedness of the Fourier transform and of the operator $T$ imply that relations \e{eq:invpr} and \e{eq:invprA} remain true for all $ u \in L_2({\Bbb R}_+)$.
 \end{remark}


\section{The perturbation determinant and the spectral shift function}
 

 \medskip
 
{\bf 6.1.}
Mathematical theory of the spectral shift function was constructed by M.   Kre\u{\i}n. Let us recall here briefly its basic notions (see \cite{BY2} and \cite{Ya}, for details). Suppose that self-adjoint operators $H_{0}$ and $H$ are semibounded from below and that the difference of their resolvents belongs to the trace class.

 The (generalized) perturbation determinant $D(z)$ for the pair $H_{0}$, $H$ is defined by the equation
\begin{equation}
 \tr (R(z)-R_{0}(z))=- D'(z) D(z)^{-1}, \q z\not\in \sigma(H_{0})\cup  \sigma(H ),
 \label{eq:D}\end{equation}
 which fixes $D(z)$ up to a constant factor. We set
  \begin{equation}
 \ln D(z)= \int_{z_{0}}^z\tr (R_{0}(z')-R (z'))dz' 
 \label{eq:4tr6}\end{equation}
 where $z_{0}$ is some real point lying below 
\[
\nu=\min\{\inf\sigma ( H_{0}), \inf \sigma (H )\} 
\]
  and the integral is taken over some contour in the upper (lower) half-plane if $\Im z>0$ ($\Im z < 0$). Then the function $\ln D(z)$ is determined up to a real constant (and hence  the function $D(z)$ is determined  up to a constant positive factor)  which is inessential. Clearly, $\arg D(z)=0$ for $ z=\bar{z}< \nu $ and $D(\bar{z})= \overline{D(z)}$.
 
  The limit of  $\arg D(\lambda+i \varepsilon)$ exists for  a.e. $ \lambda\in {\Bbb R}$, and the  spectral shift function $\xi(\lambda)$ for the pair $H_{0}$, $H$  is defined by the formula 
     \begin{equation}
\xi(\lambda) =\pi^{-1}\arg D(\lambda+i 0) . 
 \label{eq:4tr7}\end{equation}
 The function $\xi(\lambda)$ assumes constant integral values on component intervals of the set of common regular points of the operators $H_{0}$ and $H$,
   $\xi(\lambda)=0$ if $\lambda<\nu$
   and
     \begin{equation}
\xi(\lambda+0)-\xi(\lambda-0) =n_{0}-n 
 \label{eq:jump}\end{equation}
 if $\lambda $ is
 an isolated eigenvalue of multiplicity $n_{0}$ of  the operator $H_{0}$   and of multiplicity $n$ of  the operator $H$.   The spectral shift function satisfies the condition
 \[
 \int_{-\infty}^\infty |\xi(\lambda)| (1+\lambda^2)^{-1} d\lambda< \infty.
 \]
The trace formula
\[
  \tr (\varphi(H)-\varphi(H_{0}))=  \int_{-\infty}^\infty \xi(\lambda)\varphi'(\lambda) d\lambda
    \]
holds  for functions $\varphi\in C^2({\Bbb R})$ such that
  \[
  (\lambda^2 \varphi' (\lambda))'=O(\lambda^{-1-\varepsilon}), \q \lambda\to+\infty,
  \]
  for some $\varepsilon>0$. In particular, the trace formula   is true for the function $\varphi(\lambda)=(\lambda-z)^{-1}$ when $\varphi(H)=R(z)$.
  
For the pair $H_{0}$,  $H$,   wave operators \e{eq:WO} exist,  are complete and the relation between the scattering matrix $S (\lambda)$ and the spectral shift function is given by the 
Birman-Kre\u{\i}n formula  \cite{BK}
 \begin{equation}
  \det S (\lambda)=e^{-2 \pi  i \xi (\lambda)}. 
\label{eq:BK1}\end{equation}

 \medskip
 
{\bf 6.2.}
Now we are in a position to construct explicitly the perturbation determinant and the spectral shift function for the pair $H_{0}$, $H$ considered in the previous section.
According to equalities \e{eq:4ord1} and \e{eq:exP1}  we have
\begin{equation}
 \tr (R(z)-R_{0}(z))=\tr P(z)-\tr P_{0}(z).
 \label{eq:4tr}\end{equation}
 The kernel $P(x,y;z)$ of the operator $P(z)$ is determined by formulas \e{eq:4ord2} and \e{eq:4ord3} so that
 \[
 P(x,x;z)= p_{11}(\zeta) e^{2i\zeta x}+(p_{12}(\zeta)+p_{21}(\zeta)) e^{(-1+i)\zeta x} 
 +p_{22}(\zeta) e^{-2\zeta x}.
 \]
 It follows that
 \[
\tr P(z)=\int_{0}^\infty P(x,x;z) dx= (2\zeta)^{-1}( i p_{11}(\zeta)  + (1+i)(p_{12}(\zeta)+p_{21}(\zeta))   +p_{22}(\zeta) ).
\]

 Using formulas \e{eq:P} and \e{eq:4ord9} for the coefficients $p_{jl}$, we find that
  \[
\tr P(z)=  (4\zeta^4 \Omega(\zeta))^{-1}\big(\alpha_{0}+2 i \Re\alpha \zeta^2 
 + 2 (1+i) \alpha_{1} \zeta^3+ 3 \zeta^4\big).
 \]
  Similarly, it follows from  \e{eq:exP1a} that
  \begin{equation}
\tr P_{0}(z)=(4z)^{-1}. 
 \label{eq:4tr2B}\end{equation}
 Therefore equality \e{eq:4tr} yields the following result.

  \begin{proposition}\label{xi} 
  Suppose  that boundary conditions  \e{eq:4ord}  are satisfied. Then 
  \begin{equation}
 \tr (R(z)-R_{0}(z))= -4^{-1} \zeta^{-3} \Omega(\zeta)^{-1} \Omega'(\zeta) 
  \label{eq:4tr3}\end{equation}
  where the function $\Omega$ is defined by formula \e{eq:4ord8}.
  \end{proposition}
  
  Let us now consider  boundary conditions \e{eq:par2} and \e{eq:par1}. Remark that formally conditions \e{eq:par2} can be obtained from \e{eq:4ord} if we set $\alpha_{1}=N$, replace $\alpha$ by $-\alpha N$, replace  $\alpha_{2}$ by $|\alpha|^2 N + \alpha_{2}$ and pass to the limit $N\to\infty$. We plug these expressions into \e{eq:4ord8} and observe that the limit of the right-hand side of \e{eq:4tr3} is determined only by the terms containing the factor $N$. Thus, the 
   equation \e{eq:4tr3} is true if we set
    \begin{equation}
   \Omega(\zeta)=\alpha_{2}+ (1-i) |\alpha|^2 \zeta-2i \Re \alpha \, \zeta^2 -(1+i) \zeta^3.
     \label{eq:O1}\end{equation}
   Similarly, in case \e{eq:par1} we set $\alpha=0$, $\alpha_{2}=N$ and take the limit $N\to\infty$. This yields
    \begin{equation}
   \Omega(\zeta)=\alpha_1+ (1-i)   \zeta .
   \label{eq:O2}\end{equation}
   All results obtained in Section~5 for  boundary conditions \e{eq:4ord} can easily be carried over to cases  \e{eq:par2} and \e{eq:par1}. In particular, formula \e{eq:4ord17}  remains true.

Comparing equations    \e{eq:D} and \e{eq:4tr3}, we obtain

  \begin{proposition}\label{xi1}
  The perturbation determinant for the pair $H_{0}$, $H$ is given by the equality
  \begin{equation}
D(z)=  \Omega(\zeta), \q z=\zeta^4.
 \label{eq:Omega}\end{equation} 
    \end{proposition}

 Since all functions \e{eq:4ord8}, \e{eq:O1} and \e{eq:O2} satisfy the condition $\Omega(e^{\pi i/ 4} k)> 0$ for large   $k>0$, we can set $\arg \Omega(e^{\pi i/ 4} k)=0$ for such $k$. Thus, $  \Omega(\zeta)$ determines   by formulas  \e{eq:4tr7}, \e{eq:Omega}   the spectral shift function  $\xi(\lambda)$ for the pair $H_{0}$, $H$.    Clearly, $\xi(\lambda)=0$ below the lowest eigenvalue of     $H$.
 Proposition~\ref{negX} implies that the spectral shift function  has a jump $-1$
 at a simple eigenvalue of  $H$ and a jump $-2$ at an eigenvalue of multiplicity $2$. This is of course consistent with general formula \e{eq:jump}.
In view of definition   \e{eq:4tr7},  relation \e{eq:4ord17} reduces to formula \e{eq:BK1}.

 Let us write formulas \e{eq:4ord8}, \e{eq:O1} and \e{eq:O2} in a unified way as
  \begin{equation}
D(z)=\sum_{j=0}^4 \omega_{j}\zeta^j, \q z=\zeta^4.
 \label{eq:Ome1}\end{equation}
 We put
 \begin{equation}
\gamma_{0}=\min_{j}\{j: \omega_{j}\neq 0\}, \q \gamma_{1}=\max_{j}\{j: \omega_{j}\neq 0\}.
 \label{eq:Ome2}\end{equation}

 Results about the behavior of $\xi(\lambda)$  for $\lambda>0$  are collected in
 the following assertion.
 
  \begin{proposition}\label{4tr} 
  For $\lambda>0$, the spectral shift function for the pair $H_{0}$, $H$ is infinitely differentiable away from  a positive eigenvalue $\lambda_{0}$ of the operator    $H$ $($if it exists$)$. The limits of $\xi(\lambda)$ as $\lambda\to\lambda_{0}\pm 0$ exist and 
    \begin{equation}
 \xi (\lambda_{0}+ 0)- \xi (\lambda_{0} - 0)=-1.
 \label{eq:pos}\end{equation}
  If an interval $(\lambda_{1}, \lambda_{2})  \subset{\Bbb R}_{+}$ does not contain an eigenvalue of $H$, then 
      \begin{equation}
  \xi(\lambda_{2})- \xi(\lambda_{1}) =\pi^{-1} \int^{k_{2}}_{k_{1}}\Im (\Omega(k)^{-1}\Omega'(k))   dk, \q \lambda_{j}=k_{j}^4.
  \label{eq:4tr15}\end{equation}
  Moreover, there exists the limit $\xi (+0 ) $,
   \begin{equation}
\delta:=\xi(+0)- \xi(-0)= -  \gamma_{0}/4
  \label{eq:Ome}\end{equation}
   and
  \begin{equation}
 \lim_{\lambda\to +\infty}\xi (\lambda ) = - \gamma_{1}/4.
 \label{eq:pos1}\end{equation}
\end{proposition}
 
  \begin{pf} 
  According to Proposition~\ref{4ord1q} the perturbation determinant $D(\lambda+ i0)=0$ if and only if $\lambda=\lambda_{0}$ is an eigenvalue  of the operator    $H$.
  Since $D(\lambda+ i0)$ is a  $C^\infty$-function for $\lambda>0$, its argument is also   a $C^\infty$-function away from the point $\lambda_{0}$.     According to Proposition~\ref{4ord1qq},  $\lambda_{0}$ is a simple zero of the function $D(z)$ so that
  \[
 D(z)= d ( z -\lambda_{0})  +  O( | z -\lambda_{0} |^2),
  \q z \to \lambda_{0},  \q \Im z\geq 0, 
  \]
 for some $d \neq 0$. Therefore the function $\arg D(\lambda+i 0)$ has  finite limits   as $\lambda \to\lambda_{0}\pm 0$ and the variation of $\arg D(z)$ as $z$ passes from $ \lambda_{0}-\varepsilon$ to $ \lambda_{0}+\varepsilon$ in the clockwise direction over a semi-circle  
  $C_{\varepsilon}^+(\lambda_{0})= \{|z-\lambda_{0}| =\varepsilon, \Im z\geq 0\}$ equals
  \begin{equation}
 \var_{C_{\varepsilon}^+ (\lambda_{0})}\arg D(z)= -\pi+o(1) 
 \label{eq:V11}\end{equation}
  as $\varepsilon\to 0$. This proves formula \e{eq:pos}. Relation \e{eq:4tr15} is a  direct consequence of definition \e{eq:4tr7} and relation  \e{eq:Omega}.
 
The existence of $\xi(+0)$ and relation \e{eq:Ome} follow from formula  \e{eq:Ome1} and the equality
 \begin{equation}
 \var_{C_{\varepsilon}^+ (0)}\arg D(z)= -\pi\gamma_{0}/4+o(1) 
 \label{eq:V2}\end{equation}
 as $\varepsilon\to 0$. Similarly, relation \e{eq:pos1} follows from formula  \e{eq:Ome1} and the equality
 \begin{equation}
 \var_{C_{R}^+ (0)}\arg D(z)= -\pi\gamma_{1}/4 +o(1) 
 \label{eq:V3}\end{equation}
 as $R \to \infty$.
     \end{pf}

   If  boundary conditions  \e{eq:4ord}  are satisfied, then $\gamma_{1}=4$, but $\gamma_0$ might equal (see \e{eq:4ord8}) any number between $0$ and $4$. 
      In particular, if $\alpha=\alpha_{1}=\alpha_{2}=0$, then 
       $D(z)=-z $  and $\xi(\lambda)=0$ for $\lambda<0$ and  $\xi(\lambda)= -1$  for $\lambda>0$.


\section{Zero-energy resonances}
 
 
 Here we analyse resolvent singularities and the behavior of the spectral shift function at the bottom of the continuous spectrum. The notion  of a zero-energy resonance is introduced and discussed.
 
  \medskip
 
{\bf 7.1.}
Resolvents of all operators considered in Section~5 admit asymptotic expansions (for fixed $x$ and $y$) in powers of $\zeta$ as $|\zeta|\to 0$. Let us find their singular parts. Below the symbol $``\simeq"$ means an equality valid up to regular terms.

According to 
   formula \e{eq:4ordr} 
\begin{equation}
R_{00}(x,y;z) = 2 ^{-1/2}   xy (-z)^{-1/4}+ 
12 ^{-1} ( | x-y|^3-(x+y)^3) +O (|z|^{1/4}), 
  \label{eq:ju4}\end{equation}
    according to    formula \e{eq:exP1a}
      \begin{equation}
P_{0}(x,y;z) = - 2 ^{-1/2}   xy  (-z)^{-1/4}+  2^{-1}  xy  (x + y ) +O (|z|^{1/4})
  \label{eq:ju4b}\end{equation} 
   and according to    formula \e{eq:exP2}
    \begin{equation}
P_1(x,y;z)\simeq  2 ^{-1/2} (-z)^{-3/4}- 2 ^{-3/2}  (x+y)^2 (-z)^{-1/4} .  
   \label{eq:ju4c}\end{equation}
  
  It is slightly more difficult to find the singular part of kernel $P (x,y;z)$. It follows from equations   \e{eq:4ord2},  \e{eq:4ord3} and \e{eq:P} that 
  \begin{equation}
P (x,y;z) =  2^{-3/2} e^{-\pi i/4} \zeta^{-3} \Omega (\zeta)^{-1} \sum_{n=0}^\infty
(n!)^{-1} L_{n}(x,y;\zeta) \zeta^n 
  \label{eq:P1}\end{equation}
  where
    \begin{eqnarray*}
L_{n}(x,y;\zeta)& =& {\bf p}_{11}(\zeta) i^n (x+y)^n + {\bf p}_{12}(\zeta)  (-x +iy )^n 
\\
&+ &{\bf p}_{21}(\zeta)   (i x - y)^n + {\bf p}_{22}(\zeta) (-1)^n (x+y)^n
\end{eqnarray*}
   and the coefficients $ {\bf p}_{jl}(\zeta)$, $j,l=1,2$, are defined by formulas  \e{eq:4ord9}.  In particular, we have 
 \begin{equation}
L_{0}(x,y;\zeta) =  2 (1+i) \alpha_{1}\zeta ^3 + 2\zeta ^4,
  \label{eq:L0}\end{equation}  
   \begin{equation}
L_{1}(x,y;\zeta) = - 2 (1+i)( \alpha x +\bar{\alpha}y)\zeta ^2+ 2(-1+i) (x+y)\zeta ^4 
  \label{eq:L1}\end{equation} 
  and
   \begin{equation}
L_{2}(x,y;\zeta) =  -4i \alpha_{0}xy + 2 (\alpha x^2 +4 \Re\alpha x y+\bar{\alpha}y^2)\zeta ^2 + 2 (i-1) \alpha_{1} (x^2+y^2)\zeta ^3- 4i xy \zeta ^4.
  \label{eq:L2}\end{equation}  
  
  Let first $\alpha_{0}\neq 0$. Then $\Omega(0)=\alpha_{0}\neq 0$ and the only singular term in \e{eq:P1} comes from the first term in the right-hand side of  \e{eq:L2} which yields
     \begin{equation}
P(x,y;z)\simeq - 2 ^{-1/2}  xy (-z)^{-1/4}  ,\q \alpha_{0}\neq 0.
  \label{eq:ju4a}\end{equation}
  If $\alpha_{0}=0$ but $\alpha_2 \neq 0$, then ${\bf p}_{jl}(\zeta )=O (| \zeta |^2)$, $\Omega(\zeta)= (1-i) \alpha_{2}\zeta (1+ O (| \zeta |)$ as $| \zeta | \to 0$ and  singular terms in \e{eq:P1} come from the first terms in the right-hand sides of \e{eq:L0} and \e{eq:L1}. Thus, we have
     \begin{equation}
P(x,y;z)\simeq 2 ^{-3/2} \alpha_2^{-1} (2\alpha_{1}-\alpha x-\bar{\alpha}y) (-z)^{-1/4}   ,\q \alpha_{0}= 0, \q \alpha_2\neq 0.
  \label{eq:P3}\end{equation}
  If $\alpha_{0}= \alpha_2 = 0$ (and hence $\alpha = 0$) but $\alpha_{1}\neq 0$, then ${\bf p}_{jl}(\zeta )=O (| \zeta |^2)$,  $\Omega(\zeta)= -(1+ i) \alpha_{1}\zeta^3 (1+ O (| \zeta |)$ as $| \zeta | \to 0$ and  singular terms in \e{eq:P1} come from \e{eq:L0} and the second term in the right-hand side of  \e{eq:L1}. Thus, in this case we have
     \begin{equation}
P(x,y;z)\simeq 2 ^{-1/2}  (-z)^{-3/4} -(\alpha_{1}^{-1} 2^{-1/2}(x+y) +2^{-3/2}(x^2+y^2)) (-z)^{-1/4}  .
  \label{eq:P4}\end{equation}
  Let finally $\alpha=\alpha_{1}= \alpha_2 = 0$. Then  it follows from formula \e{eq:000} that
      \begin{equation}
P(x,y;z)\simeq 2 ^{-1/2}  (-z)^{-3/4} - 2^{-1}(x+y) (-z)^{-1/2} +2^{-1/2} x y  (-z)^{-1/4}  .
  \label{eq:P5}\end{equation}

 \medskip

{\bf 7.2.}
  Here we discuss different types of zero-energy resonances. 
  Let us distinguish several cases.
   
   $1^0$ The operator $H_{0}$  does not have zero-energy resonances. Since   singular terms in \e{eq:ju4} and \e{eq:ju4b} are compensated in
   \e{eq:exP1}, the resolvent kernel $R_{0}(x,y;z)$   is a continuous function as $z\to 0$ and
   \[
R_{0}(x,y;0)=
12 ^{-1} ( | x-y|^3-(x+y)^3) +  2^{-1} xy (x+ y).
  \]
  Linear functions (except zero) do not satisfy boundary conditions \e{eq:4ordfr}.
      All boundary conditions \e{eq:4ord}, \e{eq:par2}   and \e{eq:par1}  change the domain of quadratic form of the  operator $H_{0}$,  and hence the corresponding operators $H$ cannot be considered as small perturbations of $H_{0}$. 
      
         $2^0$ The operator $H_{00}$   has a zero-energy resonance. Formula \e{eq:ju4}    shows that its resolvent kernel has the  singularity    $2 ^{-1/2}   xy (-z)^{-1/4}$ as $z\to 0$.   The function  $u(x)=x$ 
satisfies  boundary conditions  \e{eq:4ordfr}. According to equality \e{eq:O2} where $\alpha_{1} = 0$
the perturbation determinant  for the pair $H_{0}$, $H_{00}$ equals  $D_{00}(z)=(-z)^{1/4}$,  and the spectral shift function      equals $0$ for $\lambda<0$   and $-1/4$  for $\lambda>0$. 
In view of formula \e{eq:ju4}  it is natural to say that the operator $H_{00}$ has  a quarter-bound state at zero energy.   The operator $H_{00}$ belongs to family \e{eq:par1}  for $\alpha_{1} = 0$. One negative eigenvalue appears for an arbitrary   $\alpha_{1}<0$.

  Next we consider the operator $H$ with boundary conditions  \e{eq:4ord}.  According to formula \e{eq:4ord1} the singular part of its resolvent equals the sum of singular parts of \e{eq:ju4} and     of the kernel $P(x,y;z)$. The jump $\delta $ of the spectral shift function $\xi(\lambda) $ at the point $\lambda=0$ is determined by formula \e{eq:Ome}.
  
    $3^0$ If $\alpha_{0}\neq 0$, then the operator $H$ does not have zero-energy resonances. Since
the singularities in \e{eq:ju4} and \e{eq:ju4a} are compensated,  the kernel $R (x,y;z)$ is a continuous function as $z\to 0$.   Linear functions $u(x)=Ax+B$   satisfy boundary conditions  \e{eq:4ord} only for $A=B  = 0$.
   Since $\Omega(0)=\alpha_{0}\neq 0$, the spectral shift function is continuous at the point  $\lambda=0$. According to Proposition~\ref{neg}   new negative eigenvalues of the operator $H$ cannot appear under small perturbations of the coefficients $\alpha$, $\alpha_1$ and $\alpha_{2}$.

    $4^0$  Let $\alpha_{0}=0$ but $\alpha_2\neq 0$.   Then
    formulas \e{eq:ju4}, \e{eq:P3} show that the resolvent kernel has a singularity at $z= 0$ of the same order  $-1/4$  as $R_{00}(z)$.  The   linear function $u(x)=\alpha_{2} x -\bar{\alpha}$  satisfies boundary conditions  \e{eq:4ord}.  Now $\omega_{0}=0$ but 
     $\omega_{1} = (1-i) \alpha_{2} \neq 0$ in \e{eq:Ome1}   and hence $\gamma_{0}=1$.  
  It follows from \e{eq:Ome} that   $\delta =-1/4 $.  Thus,  the operator $H$ has a zero-energy resonance of the same ``strength" (that is $1/4$-bound state) as the operator $H_{00}$.

$5^0$ Let $\alpha =\alpha_2=0$ but $\alpha_1\neq 0$.
 Then    formula \e{eq:P4} shows that its resolvent kernel has the singularity    $2 ^{-1/2}  (-z)^{-3/4}$ as $|z |\to 0$.  The equation $u^{(4)}(x)=0$ has  a solution $u(x)=1$  satisfying boundary condition  \e{eq:4ord}. 
   Now $\omega_{0}=\omega_{1}=\omega_{2}=0$ but $\omega_{3} = -
(1+i) \alpha_{1} \neq 0$    and hence $\gamma_{0}=3$.  
  It follows from \e{eq:Ome} that   $\delta =-3/4 $. 
To comply with  these results,
we say that the operator $H$ has a  $3/4$-bound state at energy zero. 

In the cases  $4^0$ and $5^0$, we have that $\det A=0$ but $A\neq 0$ so that the matrix $A$ has exactly one zero eigenvalue. Therefore the operator $H$ has an additional negative eigenvalue for an arbitrary negative perturbation of $A$.

$6^0$  Let $\alpha =\alpha_1=\alpha_2=0$. Then both functions $u(x)=1$ and $u(x)=x$ satisfy  boundary conditions  \e{eq:4ord}.  It follows from formula $\Omega(\zeta)=-\zeta^4$ that $\delta=-1$.  The operator $H$ has both    $3/4$- and $1/4$-bound states at energy zero. This is consistent   with formula $\delta=-1$ as well as with the fact that the operator $H$ has two negative eigenvalues for an arbitrary matrix $A<0$. According to  \e{eq:P5} the resolvent kernel contains more singularities than in the previous cases.  

Families  \e{eq:par2} and  \e{eq:par1} can be considered in a similar but simpler way.
We discuss only the operator $H_{1}$.

$7^0$   Let boundary conditions  \e{eq:ex2} be satisfied.
Formula \e{eq:ju4c}   shows that the singularity $2 ^{-1/2}  (-z)^{-3/4}$ of the resolvent kernel is the same as in case $5^0$.        The function $u(x)=1$ satisfies     conditions \e{eq:ex2}. According to equality \e{eq:O1} where $\alpha=\alpha_{2} = 0$ 
the perturbation determinant  for the pair $H_{0}$, $H_1$ equals  $D_1(z)=(-z)^{3/4}$,  and the spectral shift function      equals $0$ for $\lambda<0$   and $-3/4$  for $\lambda>0$. 
Thus,  the operator $H$ has a  $3/4$-bound state at energy zero. As follows from  \e{eq:QFBC1}, the operator $H$ corresponding to boundary conditions \e{eq:par2} where $\alpha=0$ has a negative eigenvalue for all $\alpha_{2}<0$.

 \medskip

{\bf 7.3.}
We finish with an analogue of the Levinson theorem.      Consider the closed contour which consists of a small circle $C_{\varepsilon}=\{|z|=\varepsilon\}$, a big circle $C_{R}= \{|z|=R\}$ and two intervals $(\varepsilon,R)$ lying on the upper and lower edges of the cut along $[0,\infty)$. Moreover, if the operator $H$ has a positive eigenvalue $\lambda_{0}$ we go over it by small semi-circles $C_{\varepsilon}^\pm(\lambda_{0})$ where $|z-\lambda_{0}|=\varepsilon $ and  $\pm \Im z \geq 0$ (lying in the upper and lower half-planes).  Let us pass this contour in the positive direction and apply the argument principle to the function $  D(z)$. Taking into account the direction of motion and the identity $\overline{D(z)}= D(\bar{z})$, we obtain that
 \begin{equation}
\var_{C_{\varepsilon}^+}  \arg D(z)-\var_{C_{R}^+} \arg D(z)+ \arg D(R+i0) -\arg D(\varepsilon+i0)=  \pi N_{-}
  \label{eq:Lev}\end{equation}
  where $N_{-}$ is  the   number of negative eigenvalues   of the operator $H$.  Moreover, if $H$ has a  positive eigenvalue, than the term
 \begin{equation}
  \arg D(\lambda_{0}-\varepsilon+i 0) -   \arg D(\lambda_{0}+\varepsilon+i0) +
  \var_{C_{\varepsilon}^+(\lambda_{0})} \arg D(z)
   \label{eq:Lev1}\end{equation}
  should be added to the left-hand side. 
  
  Let us choose $\arg s(\lambda)$ as a continuous function of $\lambda >0$. According to formulas \e{eq:4ord17}, \e{eq:4tr7} we have that
 \begin{eqnarray*}
  \arg D(R+i0)- \arg D(\lambda_{0}+\varepsilon+i0)&=& -2^{-1}( \arg s(R)- \arg s(\lambda_{0}+\varepsilon)),
   \\     
   \arg D(\lambda_{0}-\varepsilon+i0) -  \arg D(\varepsilon+i0)&=& -2^{-1}( \arg s(\lambda_{0}-\varepsilon)- \arg s( \varepsilon)).
 \end{eqnarray*}
     Taking the limits $\varepsilon\to 0$, $R\to\infty$ and using relations   \e{eq:V11} -- \e{eq:V3}, we finally obtain that
   \begin{equation}
 \arg s ( {\infty}) -  \arg s (+0) =-2\pi N+\pi (\gamma_{1}-\gamma_{0})/2
 \label{eq:Lev2}\end{equation}
where $N$ is the total number of eigenvalues (including eventually a positive eigenvalue) of the operator $H$ and the numbers $\gamma_{0}, \gamma_{1}$ are defined by formula  \e{eq:Ome2}. For example, for boundary conditions \e{eq:4ord}, $ \gamma_{1}=4$ and the number $\gamma_{0}$ has been computed in the previous subsection. We emphasize that an  isolated eigenvalue and an eigenvalue embedded in the continuous spectrum give the same contributions to the Levinson formula  \e{eq:Lev2}.


\begin{thebibliography}{99}
 
 
\bibitem {Ag}S. Agmon, {\em Spectral properties of Schr\"odinger operators and scattering  
 theory}, Ann. Scuola Norm. Sup. Pisa  {\bf 2} no. 4 (1975),  151-218. 
 
 \bibitem {BDT} R. Beals, P. Deift and C. Tomei, {\em Direct and inverse scattering on the line}, Math. surveys and monographs, N 28, Amer. Math. Soc., Providence, R. I., 1988.
 
 \bibitem {Bia}M. Sh. Birman, {\em On the spectrum of singular boundary-value problems}, Matem. sb. {\bf 55}, no. 2 (1961), 125-174 (Russian);
English transl., Eleven Papers on Analysis, Amer. Math. Soc. Transl. (2), vol. 53, Amer. Math. Soc., Providence, RI, 1966, 23-60. 

 
   \bibitem {B}M. Sh. Birman, {\em Scattering problems for differential operators with constant
coefficients}, Funct. Anal. Appl. {\bf 3} no. 3 (1969),  167-180.

  \bibitem {BK}M. Sh. Birman and M. G. Kre\u{\i}n, {\em On the theory of wave operators and
scattering operators},
 Soviet Math. Dokl. {\bf 3}  (1962), 740-744.
 
  \bibitem {BY4}M. Sh. Birman and D. R. Yafaev, {\em A general  scheme in the stationary theory of 
scattering}, Problemy Math. Phys. {\bf 12}  (1987), 89-117. English trasl., Amer. Math. Soc. Transl.
(Ser.2) {\bf 157}  (1993), 87-112.

\bibitem {BY2}M. Sh. Birman and D. R. Yafaev, {\em The spectral shift function. The papers of M. G.
Kre\u{\i}n and their further development},  St. Petesburg Math. J. {\bf 4} no. 5 (1993), 833-870.

 
\bibitem {CL} E. A. Coddington and N. Levinson, {\it  Theory of ordinary differential equations},
McGraw-Hill, New york, 1955.

 \bibitem {Hof} K. Hoffman, {\it Banach spaces of analytic functions}, Prentice-Hall, Inc., Englewood Cliffs, N. Y., 1962.

\bibitem {Ku}S. T. Kuroda, {\em Scattering theory for differential operators}, J. Math. Soc. Japan {\bf 25} no. 1,2 (1973), 75-104,  222-234. 

 \bibitem {Kuroda} S. T. Kuroda, {\it An introduction to scattering theory},
Lect. Notes Series No. 51, Aarhus University, 1978.



\bibitem {Pow} S. G. Power, {\em Hankel operators on Hilbert space},  Research Notes in Math., Vol 64, Pitman, Boston, 1982.
 
\bibitem {RS} M. Reed and B. Simon, {\it Methods of Modern Mathematical Physics} IV, Academic
Press, 1978.


\bibitem {F} L. D. Faddeev,  {\em Properties of the $S$-matrix of the one-dimensional Schr\"odinger equation}, Amer. Math. Soc. Transl. (Ser.2)  {\bf 65}  (1967), 139-166.

\bibitem {Y15}D. R. Yafaev, {\em The virtual level of the Schr\"odinger equation}, Zap. Nau\v{c}hn.
Sem. LOMI {\bf 51} (1975), 203-216.  English transl., J. Sov. Math. {\bf 11} (1979), 501-510.

 
 \bibitem{Ya} D. R. Yafaev, {\em Mathematical scattering theory}, Amer. Math. Soc.,   Providence,  Rhode Island, 1992. 
 
  \bibitem {Y1} D. R. Yafaev, {\em On a zero-range interaction of a quantum particle with the vacuum}, 
J. Phys. A, {\bf 25} (1992), 963-978.
 
 \bibitem {LNM}D. R. Yafaev, {\it Scattering theory: some old and new problems}, Springer Lecture Notes Math. {\bf 1735}, 2000. 
  

     
    \end{thebibliography}
\end{document}